\documentclass[12pt]{article}
\usepackage{a4,pstricks,pst-node,pst-coil,supertabular,enumerate}

\catcode`@=11
\@addtoreset{equation}{section}

\font\fiv@bb  = msbm5
\font\six@bb  = msbm6
\font\sev@bb  = msbm7   
\font\egt@bb  = msbm8
\font\nin@bb  = msbm9   
\font\ten@bb  = msbm10
\font\elv@bb  = msbm10 scaled 1096   
\font\twl@bb  = msbm10 scaled 1200

\def\tf@bb{\expandafter\csname\ifcase\@ptsize ten\or elv\or twl\fi @bb\endcsname}
\def\sf@bb{\expandafter\csname\ifcase\@ptsize sev\or egt\or egt\fi @bb\endcsname}
\def\ssf@bb{\expandafter\csname\ifcase\@ptsize fiv\or six\or sev\fi @bb\endcsname}

\newfam\@bbfam \def\bb{\fam\@bbfam \tf@bb}
\textfont\@bbfam=\tf@bb \scriptfont\@bbfam=\sf@bb \scriptscriptfont\@bbfam=\ssf@bb

\gdef\CC{{\bb C}}
\gdef\HH{{\bb H}}
\gdef\NN{{\bb N}}

\gdef\QQ{{\bb Q}}
\gdef\RR{{\bb R}}

\gdef\ZZ{{\bb Z}}

\catcode`@=12

\font\eufmtwl=eufm10 scaled 1200
\newfam\eufmfam
\textfont\eufmfam=\eufmtwl

\newfont{\goth}{eufm10 scaled 1200}
\newfont{\typewriter}{cmtt10 scaled 1200}

\newcommand{\nc}{\newcommand}
\nc{\PsoM}{\mbox{{$P_{\mbox{\scriptsize SO}}(M)$}}}
\nc{\PsonM}{\mbox{{$P_{\mbox{\scriptsize SO(n)}}(M)$}}}
\nc{\PsoG}{\mbox{{$P_{\mbox{\scriptsize SO}}(G)$}}}
\nc{\PspinM}{\mbox{{$P_{\mbox{\scriptsize Spin}}(M)$}}}
\nc{\PspinG}{\mbox{{$P_{\mbox{\scriptsize Spin}}(G)$}}}
\nc{\Pspineps}{\mbox{{$P_{\mbox{\scriptsize Spin,}\epsilon}$}}}
\def\psp{{P_{\mbox{\scriptsize Spin}}}}
\def\pso{{P_{\mbox{\scriptsize SO}}}}
\nc{\timesspinn}{\mbox{{$\times_{\mbox{\scriptsize Spin(n)}}$}}}
\nc{\SOn}{\mbox{{$\mbox{SO}(n)$}}}
\nc{\Spinn}{\mbox{{$\mbox{Spin}(n)$}}}
\def\SO{\mathop{{\rm SO}}}
\def\Spin{\mathop{{\rm Spin}}}
\def\Cl{\mathop{{\rm Cl}}}
\def\GaM#1#2{{{\ti\Ga}_{#1}^{#2}}}
\def\GaN#1#2{{{   \Ga}_{#1}^{#2}}}
\def\al{{\alpha}}
\def\om{{\omega}}
\def\na{{\nabla}}
\def\la{{\lambda}}
\def\ka{{\kappa}}
\def\si{{\sigma}}
\def\Si{{\Sigma}}
\def\ga{{\gamma}}
\def\ep{{\varepsilon}}
\def\Ga{{\Gamma}}
\def\Tau{{\cal T}}
\def\pa{{\partial}}
\def\th{{\vartheta}}
\def\ph{{\varphi}}
\def\Ph{{\Phi}}
\def\de{{\delta}}

\def\el{{\ell}}
\def\LL{{\cal L}}

\def\gg{\mbox{\goth g}}
\def\zz{\mbox{\goth z}}

\def\cZ{{\cal Z}}
\def\ker{\mathop{{\rm ker}}}
\def\image{\mathop{{\rm image}}}
\def\grad{\mathop{{\rm grad}}}
\def\span{\mathop{{\rm span}}}

\def\spec{\mathop{{\rm spec}}}
\def\dim{\mathop{{\rm dim}}}
\def\vol{\mathop{{\rm vol}}}
\def\dvol{\mathop{{\rm dvol}}}

\def\id{\mathop{{\rm id}}}
\def\Id{\mathop{{\rm Id}}}
\def\End{\mathop{{\rm End}}}
\def\res#1#2{{#1}\lower .11ex\hbox{$|$}\lower .644ex\hbox{$\scriptstyle #2$}}
\def\stelle#1#2{\left. {#1}\right|_{#2}}
\def\ohne{-}
\def\ti{\tilde}
\def\gdw{\Longleftrightarrow}
\def\bs{\backslash}
\def\ol{\overline}
\def\Hp#1{H^{\mbox{\scriptsize proj}}_{#1}}

\def\komment#1{}
\def\proof#1{{\noindent {\bf Proof {#1} }}}
\DeclareSymbolFont{Extrasymb}{U}{msa}{m}{n}
\DeclareMathSymbol\square\mathrel{Extrasymb}{"03}
\def\qed{{\leavevmode\unskip\nobreak\hfil\penalty 50\hskip 1em%
  \hbox{}\nobreak\hfil\lower 1pt\hbox{$\square $\kern-.5pt}\parfillskip 0pt
  \finalhyphendemerits 0\par}}
\newtheorem{theorem}{\bf T{\footnotesize HEOREM}}[section]
\newtheorem{lemma}[theorem]{\bf L{\footnotesize EMMA}}
\newtheorem{corollary}[theorem]{\bf C{\footnotesize OROLLARY}}
\newtheorem{remark}[theorem]{\bf R{\footnotesize EMARK}}

\begin{document}

\psset{unit=0.4cm}


\title{The Dirac Operator on Nilmanifolds and Collapsing Circle Bundles\\
\vskip.8cm\normalsize{by}}
\author{Bernd Ammann and Christian B\"ar}
\date{November, 1997}
\maketitle

\begin{abstract}
We compute the spectrum of the Dirac operator on 3-dimensional
Heisenberg manifolds.
The behavior under collapse to the 2-torus is studied.
Depending on the spin structure either all eigenvalues tend to
$\pm\infty$ or there are eigenvalues converging to those of the
torus.
This is shown to be true in general for collapsing circle bundles 
with totally geodesic fibers.
Using the Hopf fibration we use this fact to compute the Dirac
eigenvalues on complex projective space including the multiplicities.

Finally, we show that there are 1-parameter families of Riemannian
nilmanifolds such that the Laplacian on functions and the Dirac
operator for certain spin structures have constant spectrum while 
the Laplacian on 1-forms and the Dirac operator for the other 
spin structures have nonconstant spectrum.
The marked length spectrum is also constant for these families.


{\bf Mathematics Classification:}
58G25, 58G30, 53C20, 53C25, 53C30
\end{abstract}


\setcounter{section}{-1}\section{Introduction}

The Dirac operator on a Riemannian spin manifold is a natural elliptic
first order differential operator.
If the underlying manifold is closed, then the spectrum is discrete.
In order to understand the relation of the Dirac eigenvalues and 
geometric data it is desirable to have examples for which one can
explicitly compute the spectrum.
These examples can be used to test conjectures.
For example, they show that the spectrum does not determine the
manifold, not even topologically.
There are examples of Dirac isospectral pairs of spherical space forms
with different fundamental groups \cite{b3}.

In the following table we collect those closed Riemannian spin manifolds
known to us for which the Dirac spectrum has been computed.

\begin{center}
\begin{tabular}{|c|c|c|}
\hline
$\RR^n/\Gamma$ & flat tori & \cite{fried1} \\ \hline
$S^n$ & spheres of constant curvature & \cite{sul}, \cite{b3},
\cite{traut, traut2}, \cite{camhi} \\ \hline
$S^n/\Gamma$ & spherical space forms & \cite{b3} \\ \hline
$S^{2m+1}$ & spheres with Berger metrics & \cite{hi} for $m=1$ \\
&& \cite{b2} for general $m$ \\ \hline
$S^3/\ZZ_k$ & 3-dim.~lens spaces with Berger metric & \cite{b4}
\\ \hline
$G$ & simply connected compact Lie groups & \cite{fegan} \\ \hline
$\CC P^{2m-1}$ & complex projective spaces & \cite{cfg1,cfg2},
\cite{sese} \\ \hline
$\HH P^m$ & quaternionic projective spaces & \cite{bunke} for $m=2$\\
&&\cite{milh1} for general $m$\\ \hline
$Gr_2(\RR^{2m})$ & real Grassmannians & \cite{strese2} for $m=3$ \\
&& \cite{strese1} for general $m$ \\ \hline
$Gr_{2p}(\RR^{2m})$ & real Grassmannians & \cite{seeger1,seeger2}\\ \hline
$Gr_2(\CC^{m+2})$ & complex Grassmannians & \cite{milh2} \\ \hline
$G_2/SO(4)$ & & \cite{seeger1,seeger2} \\ \hline
\end{tabular}
\end{center}

We compute the Dirac spectrum for another class of Riemannian spin
manifolds, for the 3-dimensional {\it Heisenberg manifolds}.
They are of the form $M = \Ga\backslash G$ where $G$ is the 
3-dimensional Heisenberg group and $\Ga$ is a cocompact lattice.
The metric comes from a left invariant metric on $G$.
We restrict to three dimensions mostly for simplicity.
The methods also work in higher dimensions.

After collecting some facts from spin geometry and from the
representation theory of nilpotent groups in the first two sections
we perform the computation of the Dirac eigenvalues in Section 3
(Theorem \ref{heisenspec}).
Our Heisenberg manifolds are circle bundles over the flat 2-torus.
The explicit formulas for the Dirac eigenvalues allow us to study
their behavior under collapse to the 2-torus, i.e.\ we let the
length of the fibers tend to zero.
It turns out that the behavior of the eigenvalues depends essentially
on the spin structure.
For some spin structures all eigenvalues tend to $\pm\infty$.
For the other spin structures most of the eigenvalues also tend 
to $\pm\infty$ but some converge to exactly the eigenvalues of the 
torus (Corollary \ref{heisencollaps}).

As it turns out in Section 4 this is a general fact for the Dirac
operator on the total space of a collapsing circle bundle with totally
geodesic fibers over any Riemannian manifold.
There are two types of spin structures which we call {\it projectable}
and {\it nonprojectable}.
In the nonprojectable case all eigenvalues tend to $\pm\infty$
(Theorem~\ref{collapsnoproj})
whereas in the projectable case there are also eigenvalues converging
to exactly the eigenvalues of the base (Theorem~\ref{collapstheo}).

As an example we apply this to the Hopf fibration over the complex
projective space $\CC P^m$.
Using explicit formulas for the Dirac eigenvalues on Berger spheres
\cite{b2} we see that for $m$ even (when $\CC P^m$ is nonspin) all
eigenvalues go to $\pm\infty$.
If $m$ is odd (when $\CC P^m$ is spin) the limits of the eigenvalues
give us a new computation of the Dirac spectrum of $\CC P^m$
confirming (and simplifying) the results in \cite{cfg1,cfg2,sese}.
We also obtain formulas for the multiplicities (Theorem 
\ref{cpnspektrum}) which have not been computed before.

In the last section we discuss the change of the spectrum of the 
Dirac operator on certain pairs or continuous families of 
Riemannian nilmanifolds. 
A {\it Riemannian nilmanifold} is a nilmanifold $\Ga \backslash G$ 
endowed with a Riemannian metric whose lift to $G$ is left invariant. 
Here we assume that $G$ is a simply connected nilpotent 
Lie group and that $\Ga$ is a cocompact discrete subgroup.
This generalizes Heisenberg manifolds.  
Besides the Dirac spectrum we look at several different spectral 
invariants which have already been studied by other authors: 
the spectrum of
the Laplacian on functions, the spectrum of the Laplacian on forms,
the length spectrum, i.e.\ the set of all lengths of closed geodesics, 
and the marked length spectrum, i.e.\ the set of all lengths of 
closed geodesics in each free homotopy class.

The questions behind the study of these examples are:
``How far is the geometry on a Riemannian nilmanifold determined
by these spectral invariants?'' and
``How are these invariants related to each other?''
The class of nilmanifolds provides us 
with many examples for the discussion of these questions.

Milnor gave in \cite{miln} the first example of a pair of 
Riemannian manifolds (16-dimensional tori) which are 
non\-isometric but isospectral for the Laplace operator on functions.

Gordon and Wilson \cite{gw1} constructed a continuous 
family of nonisometric Riemannian 2-step nilmanifolds with the same 
spectrum of the Laplacian on functions and forms. 
This is the first example 
of a nontrivial Laplace-isospectral deformation. 
This family and other deformations of the Gordon-Wilson type, i.e. via
almost inner automorphisms, do not only have the same spectrum for 
the Laplacian on functions and on forms, but also the same 
marked length spectrum.
In the 2-step case with a left invariant metric Eberlein 
\cite[Theorem~5.20]{eber}
showed that two 2-step nilmanifolds having the same marked length 
spectrum arise via almost inner automorphisms and therefore 
have the same spectrum of the Laplacian on functions and on forms.

This statement is false if we drop the condition ``2-step''.
Gornet \cite{g2} constructed two families of deformations of 
3-step Riemannian nilmanifolds
which are isospectral for the Laplacian on functions and 
for the marked length spectrum,
but not for the spectrum of the Laplacian on 1-forms.
We will show that these families are Dirac isospectral for the
nonprojectable spin structures while they have nonconstant Dirac
spectrum for the projectable spin structures.


\section{Review of Spin Geometry}

In this section we collect a few facts on spin manifolds and the Dirac 
operator which we will use later on. 
For details the reader may consult \cite{bgv} or \cite{lm}. 
Let $M$ be an oriented Riemannian manifold of dimension $n$. 
Let $\PsoM$ be the set of all oriented orthonormal bases of all tangent 
spaces of $M$. 
Obviously, $\PsoM$ is an $\SOn$-principal bundle. 
A {\it spin structure} on $M$ consists of a $\Spinn$-principal bundle 
$\PspinM$ over $M$ together with a twofold covering map 
$\varphi:\PspinM\rightarrow\PsoM$ such that the diagram 
$$
\begin{array}{cccl}
\PspinM \times \Spinn &\rightarrow& \PspinM & \\
& & &\searrow \\
\downarrow\varphi\times\Theta & &   \downarrow\varphi & \quad M\\
& & &\nearrow \\
\PsoM \times \SOn &\rightarrow& \PsoM &
\end{array}
$$
commutes.
 Here $\Theta:\Spinn\rightarrow\SOn$ is the standard twofold 
covering map and the horizontal arrows are given by the principal bundle 
structure. 

A manifold which admits such a spin structure is called a 
{\it spin manifold}. 
Being spin is a global condition which is equivalent to the vanishing of 
the second Stiefel-Whitney class of the tangent bundle, $w_2(M)=0$. 
The number of different spin structures, if any, is given by 
$\# H^1(M, \ZZ_2)$. 
In particular, a simply connected manifold has at most one spin structure. 

Let us from now on assume that $M$ is spin and that a spin structure on 
$M$ is fixed. 
There is a unitary representation $\Sigma_n$ of $\Spinn$ of dimension 
$2^{[n/2]}$, called the {\it spinor representation}. 
If $n$ is odd, then $\Sigma_n$ is irreducible. 
The associated hermitian vector bundle 
$\Sigma M:=\PspinM \timesspinn \Sigma_n$ is called the {\it spinor bundle}.

The Levi-Civita connection on $M$ induces a natural connection 
$\na^{\Si}$ on $\Si M$. 
It can be described as follows.
Let $e_1, \ldots, e_n$ be a local orthonormal tangent frame defined over 
$U\subset M$. 
Let $\Gamma^k_{ij} : U\rightarrow\RR$ be the corresponding Christoffel 
symbols, $\nabla_{e_i} e_j = \sum\limits^n_{k=1} \Gamma^k_{ij} e_k$. 
Then $(e_1, \ldots, e_n)$ is a local section of $\PsoM$. 
Let $q$ be a lift to $\PspinM$, i.e.~$\varphi\circ q = (e_1, \ldots, e_n)$. 
Then $q$ defines a trivialization of $\Sigma M$ over $U$, 
$\Sigma M|_U=U\times\Sigma_n$, with respect to 
which we have the following formula for $\nabla^{\Sigma}$:
\begin{equation}\label{nablaspinformel}
\nabla^{\Sigma}_{e_i} \sigma = \partial_{e_i} \sigma + 
\frac{1}{4} \sum\limits^n_{j, k=1} \Gamma^k_{ij} \gamma(E_j) \gamma(E_k) 
\sigma
\end{equation}

Since the spinor representation of $\Spinn$ extends to a representation 
of the Clifford algebra $\Cl(n)$ there is a well-defined 
{\it Clifford multiplication} 
$T_p M \otimes \Sigma_p M \rightarrow \Sigma_p M$, 
$X\otimes \sigma\rightarrow\gamma(X)\sigma$. 
It satisfies the relations 
$\gamma(X)\gamma(Y) + \gamma(Y)\gamma(X) + 2
\langle X, Y\rangle\, \mbox{id} =0$. 

For example in dimension $n=3$, one can choose a basis for $\Sigma_3$ 
with respect to which 
$$
\gamma(E_1) = \left(\begin{array}{cc}  
0 & i \\ i & 0 
\end{array}\right), 
\gamma(E_2) = \left(\begin{array}{cc}  
0 & -1 \\ 1 & 0 
\end{array}\right), 
\gamma(E_3) = \left(\begin{array}{cc}  
i & 0 \\ 0 & -i 
\end{array}\right),
$$
where $E_1, E_2, E_3$ denotes the standard basis of $\RR^3 \subset \Cl(3)$. 

The Dirac operator $D$ acts on the sections of $\Sigma M$. 
It is defined by $D\sigma = \sum\limits^n_{i=1} \gamma(e_i) 
\nabla^{\Sigma}_{e_i}\sigma$ where $e_1, \ldots, e_n$ is any orthonormal 
basis of the tangent space. 
If one has, in addition, a complex vector bundle $E$ over $M$ equipped
with a connection one can form the {\em twisted Dirac operator} $D^E$
acting on sections of $\Sigma M \otimes E$ by using the tensor product
connection and Clifford multiplication on the first factor.

The Dirac operator is a formally self adjoint elliptic differential 
operator of first order. 
If the manifold $M$ is closed, then $D$ has discrete real spectrum. 

Let $G$ be a simply connected Lie group with a left invariant metric. 
We regard the elements of the Lie algebra {\goth g} as left invariant 
vector fields. 
The choice of an oriented orthonormal basis of {\goth g} yields a 
trivialization of the frame bundle $\PsoG = G\times \SOn$. 
The unique spin structure can be written as $\PspinG = G\times \Spinn$, 
$\varphi= \mbox{id} \times \Theta$. 
Spinor fields are then simply maps $G\rightarrow\Sigma_n$. 

Let $\Gamma\subset G$ be a lattice. 
Spin structures of $M = \Gamma\setminus G$ correspond to homomorphisms 
$\epsilon:\Gamma\rightarrow\ZZ / 2 \ZZ = \{-1, 1\}$. 
The corresponding spin structure is given by 
$\Pspineps (\Gamma\setminus G) = G\times_{\Gamma} \Spinn$ where 
$g_0\in\Gamma$ acts on $G$ by left multiplication and on $\Spinn$ by 
multiplication with the central element $\epsilon(g_0)$. 
Spinor fields on $M$ can then be identified with $\epsilon$-equivariant 
maps $\sigma : G\rightarrow\Sigma_n$, 
i.e.~$\sigma(g_0 g) = \epsilon(g_0)\sigma(g)$ for all $g\in G$, 
$g_0 \in \Gamma$. 
Denote the corresponding Hilbert space of square integrable spinors
by $L^2(\Sigma_\ep M)$.


\section{Review of Kirillov Theory}
In this section we summarize some facts from the representation theory of 
nilpotent groups which will be of importance for our study of 
nilmanifolds. 
For details see \cite{cg}. 
Let $G$ be a simply connected nilpotent Lie group with Lie algebra 
{\goth g}. The exponential map $\mbox{exp}: \mbox{{\goth g}}\rightarrow G$ is a 
global diffeomorphism whose inverse we denote by 
$\log: G\rightarrow\mbox{{\goth g}}$. 

\noindent
{\bf Example.} Let 
$$G=\left\{\left(\left.
\begin{array}{ccc}
1 & x & z\\
0 & 1 & y\\
0 & 0 & 1
\end{array}
\right) 
\right| 
x, y, z\in\RR
\right\}
$$
be the 3-dimensional {\it Heisenberg group}. Its Lie algebra is 
$$
\mbox{{\goth g}} = 
\left\{\left(\left.
\begin{array}{ccc}
0 & x & z\\
0 & 0 & y\\
0 & 0 & 0
\end{array}
\right) 
\right| x, y, z\in\RR
\right\}, 
$$ 
the 3-dimensional {\it Heisenberg algebra}. 
Write
$$
\begin{array}{rcl}
g(x, y, z) &:=&
\left(
\begin{array}{ccc}
1 & x & z\\
0 & 1 & y\\
0 & 0 & 1
\end{array}
\right) \quad \mbox{and} \\
\\
X(x, y, z) &:=& 
\left(
\begin{array}{ccc}
0 & x & z\\
0 & 0 & y\\
0 & 0 & 0
\end{array}
\right).
\end{array}
$$
Then the exponential map is given by 
$$
\begin{array}{rcl}
\mbox{exp} (X(x,y,z)) &=& g \left(x,y,z + \frac{1}{2} x\cdot y\right)
\quad \mbox{and its inverse by} \\
\log (g(x,y,z)) &=& X \left(x,y,z - \frac{1}{2} x\cdot y\right).
\end{array}
$$
Define generators of {\goth g} by
$$
\begin{array}{rcl}
X &:=& X (1,0,0), \\
Y &:=& X (0,1,0), \\
Z &:=& X (0,0,1). \\
\end{array}
$$
Then $[X, Y] = Z$ ist the only nontrivial commutator and $Z$ spans the 
center of {\goth g}. 

Any Lie group acts on its Lie algebra via the {\it adjoint representation} 
$$
\mbox{Ad} : G \rightarrow \mbox{End} (\mbox{{\goth g}}), \, 
\mbox{Ad}_g (Y) = \frac{d}{dt} \left.\left(g\cdot \mbox{exp} (tY)\cdot
g^{-1}
\right)\right|_{t=0}
$$
and on the dual space $\mbox{{\goth g}}^*$ via the 
{\it coadjoint representation} 
$$
\mbox{Ad}^* : G\rightarrow \mbox{End} (\mbox{{\goth g}}^*), \, 
(\mbox{Ad}^*_g l) (Y) = l(\mbox{Ad}_{g^{-1}} Y), \, 
g\in G, Y\in \mbox{{\goth g}}, l\in \mbox{{\goth g}}^*.
$$
The orbits for the coadjoint representation are called 
{\it coadjoint orbits}. 
Given any $l\in \mbox{{\goth g}}^*$ there exists, in the nilpotent case, a 
maximal isotropic subspace 
$\mbox{{\goth m}} \subset \mbox{{\goth g}}$ 
for the antisymmetric bilinear form $(X, Y) \rightarrow l([X, Y])$ which 
is also a subalgebra. 
Such subalgebras are called {\it polarizing subalgebras.} 

\noindent
{\bf Example.} Let $G$ be the 3-dimensional Heisenberg group. 
Let $l\in \mbox{{\goth g}}^*$. 
We have to distinguish two cases depending on whether $l$ vanishes on the 
center or not. 

{\bf Case 1:} $l(Z)= 0$. 
$$
\begin{array}{rcl}
\left(\mbox{Ad}^*_{g (x, y, z)} l\right) 
\left( X(\xi, \eta, \zeta)\right) 
&=& l \left(\mbox{Ad}_{g (-x, -y, -z +xy)} X(\xi, \eta, \zeta)\right) \\
&=& l\left( X(\xi, \eta, \zeta + \xi y - x\eta)\right) \\
&=& l(\xi\cdot X+\eta\cdot Y+(\zeta+\xi y - x\eta)\cdot Z) \\
&=& l(\xi\cdot X + \eta \cdot Y) \\
&=& l\left(X(\xi, \eta, \zeta)\right).
\end{array}
$$
Hence $\mbox{Ad}^*_g \, l=l$ for all $g\in G$, i.e.~the orbit of $l$ 
consists of $l$ only. 
The bilinear form $(X_1, X_2) \rightarrow l ([X_1, X_2])$ vanishes in 
this case. 
Hence $\mbox{{\goth m}} = \mbox{{\goth g}}$ is the unique polarizing 
subalgebra. 

{\bf Case 2:} $l(Z)= \tau \not= 0$. 

Write $l=\alpha\cdot X^* + \beta\cdot Y^* + \tau\cdot Z^*$ where 
$X^*, Y^*, Z^* \in \mbox{{\goth g}}^*$ is the basis dual to $X,Y,Z$. 
The computation above shows 
$$
\begin{array}{rcl}
\mbox{Ad}^*_{g (x, y, z)} l 
&=& (\alpha + \tau y) \cdot X^* + (\beta - \tau x) \cdot Y^* + \tau 
\cdot Z^* \\
&=& l+\tau y\cdot X^* - \tau x \cdot Y^*.
\end{array}
$$
Hence the coadjoint orbit through $l$ is the affine hyperplane through 
$l$ spanned by $X^*$ and $Y^*$. As a polarizing subalgebra we can choose 
e.g.~$\mbox{{\goth m}} = \RR\cdot X\oplus\RR\cdot Z$ or 
$\mbox{{\goth m}} = \RR\cdot Y\oplus\RR\cdot Z$. 

For a general simply connected nilpotent Lie group $G$ pick an 
$l\in\mbox{{\goth g}}^*$. 
We will construct a continuous unitrary representation of $G$ associated 
to $l$. 
Choose a polarizing subalgebra $\mbox{{\goth m}} \subset \mbox{{\goth g}}$ 
for $l$. 
Let $M := \mbox{exp} (\mbox{{\goth m}})$ be the corresponding subgroup of 
$G$. 
By the definition 
$$
\tilde{\rho}_l (\mbox{exp} (X)) := e^{2\pi il(X)}
$$ 
we obtain a well-defined 1-dimensional unitary representation 
$\tilde{\rho}_l$ of $M$. Induction yields a continuous unitary 
representation $\rho_l$ of $G$, 
$\rho_l = \mbox{Ind} (M\uparrow G, \, \tilde{\rho}_l)$. 
One checks that $\rho_l$ does not depend (up to equivalence) upon the 
choice of {\goth m}, that $\rho_{l_1}$ is equivalent to $\rho_{l_2}$ if 
and only if $l_1$ and $l_2$ are in the same coadjoint orbit and that 
every irreducible unitary representation of $G$ is equivalent to $\rho_l$ 
for some $l$. In other words, Kirillov theory sets up a bijection between 
the set of coadjoint orbits and $\hat{G}$, the set of equivalence classes 
of irreducible unitary representations of $G$. 

\noindent
{\bf Example.} Let $G$ be the 3-dimensional Heisenberg group, let 
$l=\alpha X^* + \beta Y^* + \tau Z^* \in \mbox{{\goth g}}^*$. 
Again, we have to distinguish the two cases $\tau=0$ and 
$\tau\not= 0$. 

{\bf Case 1.} $\tau=0$. 

In this case $M=G$, hence $\rho_l = \tilde{\rho}_l$ is a 1-dimensional 
irreducible representation of $G$. 
Since coadjoint orbits of those $l$ with $\tau=0$ are just points the 
corresponding representations are parametrized by the two parameters 
$\alpha$, $\beta\in\RR$, $\rho_l =: \rho_{\alpha, \beta}$, where 

\begin{eqnarray}\label{rhoab}
\rho_{\alpha, \beta} (g(x,y,z)) 
&=& e^{2\pi i \left(\alpha X^* + \beta Y^* + \tau Z^*\right) 
(X(x,y,z - \frac{1}{2} xy))} \\ 
&=& e^{2\pi i (\alpha x + \beta y)}\nonumber
\end{eqnarray}

{\bf Case 2.} $\tau\not= 0$. 

Since the representation $\rho_l$ does not change if we replace $l$ by 
another linear form in its coadjoint orbit we may assume $l=\tau\cdot 
Z^*$. 
Choose the polarizing subalgebra 
$\mbox{{\goth m}} = \RR \cdot Y\oplus\RR\cdot Z$. 

The representation space of the induced representation $\rho_l$ is the 
space of $L^2$-sections over the homogeneous space $G/M$ in the line 
bundle associated with $\tilde{\rho}_l$, the Hilbert space 
$L^2 (G/M, G\times_{\tilde{\rho}_l} \CC)$. 
Sections in $G\times_{\tilde{\rho}_l} \CC$ can be identified with 
$\tilde{\rho}_l$-equivariant maps $f:G\rightarrow\CC$, 
i.e.~$f(gm^{-1}) = \tilde{\rho}_l (m) f(g)$ for all $m\in M$, $g\in G$. 
There is a diffeomorphism $\RR\rightarrow G/M$, 
$t\rightarrow [g(t, 0, 0)]$, where $[\cdot]$ denotes the equivalence 
class in $G/M$. 
Its inverse is simply given by $[g(x,y,z)]\rightarrow x$. 
Using this diffeomorphism we can construct out of a complex valued 
$L^2$-function $u\in L^2(\RR, \CC)$ a $\tilde{\rho}_l$-equivariant map 
$f_u : G\rightarrow\CC$ by the formula
$$
\begin{array}{rcl}
f_u(g(x,y,z)) 
&=& f_u (g(x,0,0) \cdot g(0,y,z - xy) \\
&=& \tilde{\rho}_l \left( g(0,y,z - xy)^{-1}\right) \cdot f_u (g(x,0,0)) \\
&=& \tilde{\rho}_l (g(0, -y, -z +xy) \cdot u(x) \\
&=& e^{2\pi il (X(0, -y, -z +xy))} \cdot u(x) \\
&=& e^{-2\pi i\tau (z -xy)} \cdot u(x).
\end{array}
$$
Using the isomorphism 
$L^2(\RR, \CC) \rightarrow L^2(G/M, G\times_{\tilde{\rho}_l} \CC)$, 
$u\rightarrow f_u$, we can replace 
$L^2(G/M, G\times_{\tilde{\rho}_l} \CC)$ by $L^2 (\RR, \CC)$ as the 
representation space for $\rho_l$. 

How does a given element $g(x,y,z)\in G$ act on $u\in L^2(\RR, \CC)$ via 
$\rho_l$? 
On the corresponding $f_u$ it acts by 
$$
\begin{array}{rcl}
(\rho_l(g(x,y,z))f_u) (g(\tilde{x}, \tilde{y}, \tilde{z})) 
&=& f_u (g(x,y,z)^{-1} \cdot g(\tilde{x}, \tilde{y}, \tilde{z})) \\
&=& f_u (g(-x + \tilde{x}, -y + \tilde{y}, -z + \tilde{z} + xy  - 
    x\tilde{y})) \\
&=& e^{-2\pi i \tau (-z + \tilde{z} - \tilde{x} (\tilde{y} -y))} \cdot 
    u(\tilde{x} -x).
\end{array}
$$
Putting $\tilde{y}=\tilde{z}=0$ yields the corresponding function in 
$L^2(\RR, \CC)$
$$
(\rho_l(g(x,y,z))u)(\tilde{x}) = e^{-2\pi i\tau(-z+\tilde{x}y)} 
\cdot u(\tilde{x}-x).
$$
Denote the representation $\rho_l$ considered as acting on $L^2(\RR, 
\CC)$ by $\rho_{\tau}$. 

Summarizing we see that the 3-dimensional Heisenberg group has two 
families of irreducible unitary representations, the first one 
parametrized by $\alpha, \beta\in\RR$ and acting on $\CC$: 
$$
\rho_{\alpha, \beta} (g(x,y,z)) = e^{2\pi i (\alpha x + \beta y)}.
$$
The second family is parametrized by $\tau\in\RR - \{0\}$ and acts on 
$L^2 (\RR, \CC)$:
\begin{equation}\label{rhot}
\left(\rho_{\tau}(g(x,y,z))u\right) (t) = e^{2\pi i\tau (z -ty)} u(t-x).
\end{equation}
Note that $\rho_{\tau}$ is determined by the action of the center of $G$ 
$$
\rho_{\tau} (g(0,0,z)) = e^{2\pi i\tau z} \cdot \mbox{Id}.
$$
We compute the action of the Heisenberg algebra {\goth g} on 
$L^2 (\RR, \CC)$: 
$$
\begin{array}{rcl}
(\rho_{\tau})_* (X) u (t) 
&=& \frac{d}{ds} \rho_{\tau} (\mbox{exp} s X) u(t) |_{s=0} \\
&=& \frac{d}{ds} \rho_{\tau} (g(s,0,0)) u(t) |_{s=0} \\
&=& \frac{d}{ds} e^{2\pi i\tau \cdot 0} u(t -s) |_{s=0} \\ 
&=& -u' (t).
\end{array}
$$
Similarly, 
$$
\begin{array}{rcl}
(\rho_{\tau})_* (Y) u(t) &=& -2\pi i\tau \, t\, u(t) \\
(\rho_{\tau})_* (Z) u(t) &=& 2\pi i\tau  \, u(t).
\end{array}
$$
The Hilbert space $L^2 (\RR, \CC)$ has the well-known basis of 
{\it Hermite functions} 
$$
h_k (t) = e^{\raise2pt \hbox{$t^2\over2$}} \left(\frac{d}{dt}\right)^k e^{-t^2}, \, k\in\NN_0.
$$
They satisfy the relations 
$$
\begin{array}{l}
h^{'}_k (t) = t\, h_k (t) + h_{k+1} (t), \\
h_{k+2} (t)+ 2 t\, h_{k+1} (t) + 2(k+1) h_k (t) = 0.
\end{array}
$$
We look at the $L^2$-basis $u_k (t) = h_k (\sqrt{2\pi|\tau|} t)$. 
Then the relations translate into 

\begin{eqnarray}
2\pi |\tau| t \, u_k (t) + \sqrt{2\pi|\tau|}\,
u_{k+1}(t) &=& u'_k (t), 
\label{hermite1}\\
u_{k+2} (t) + 2 \sqrt{2\pi|\tau|}\, t \, u_{k+1} (t) + 2(k+1) u_k
(t)&=&0 \label{hermite2}.
\end{eqnarray}

\noindent
This basis will be used in the next section for the computation of the
Dirac spectrum of 3-dimensional Heisenberg manifolds.


\section{Heisenberg Manifolds}
In this section let $G$ denote the 3-dimensional Heisenberg group. 
Let $r$ be a positive integer. 
Define 
$$
\Gamma_r := \left\{ g (r \cdot x,y,z) \in G \mid x,y,z \in \ZZ \right\}.
$$
Then $\Gamma_r$ is a uniformly discrete cocompact lattice in $G$.
For every lattice in $G$ there is an automorphism of $G$ mapping the 
lattice to some
$\Gamma_r$, \cite[\S 2]{gw}. 
We want to compute the Dirac spectrum of the Heisenberg manifold 
$M=\Gamma_r\bs G$. 

For positive reals $d$, $T>0$ we equip $G$ and $M$ with the left 
invariant metric for which 
$$
\begin{array}{rcl}
e_1 &:=& -d \cdot X, \\
e_2 &:=& -d \cdot Y, \\
e_3 &:=& T^{-1} \cdot Z \\
\end{array}
$$
form an orthonormal frame. 
Here $X, Y, Z$ denote the standard generators of {\goth g} (compare last 
section). 
Using the Koszul formula we compute the Christoffel symbols 
$$
\Gamma^3_{12} = \Gamma^1_{23} = \Gamma^1_{32} = 
- \Gamma^3_{21} = - \Gamma^2_{31} = - \Gamma^2_{13} = 
\frac{d^2 T}{2}, \Gamma^k_{ij}=0 \quad \mbox{ otherwise}.
$$
To emphasize the dependence of $M$ on the parameters $r$, $d$, and $T$ we 
also write $M(r,d,T)$ instead of $M$. 
Note that $d$ and $T$ are differential geometric parameters whereas $r$ 
is topological. 
The fundamental groups $\pi_1 (M(r,d,T))=\Gamma_r$ are nonisomorphic for 
different values of $r$. 
We regard spinor fields on the universal covering as maps 
$G\rightarrow\Sigma_3$. 
Choose a basis of $\Sigma_3 \cong \CC^2$ such that Clifford 
multiplication of the standard basis $E_1, E_2, E_3 \in \RR^3$ is given 
by 
$$
\gamma(E_1) = \left(\begin{array}{cc}  
0 & i \\ i & 0 
\end{array}\right), 
\gamma(E_2) = \left(\begin{array}{cc}  
0 & -1 \\ 1 & 0 
\end{array}\right), 
\gamma(E_3) = \left(\begin{array}{cc}  
i & 0 \\ 0 & -i 
\end{array}\right),
$$
compare first section.
The spinor connection on $G$ with respect to the left invariant metric 
described above is given by (\ref{nablaspinformel})
$$
\begin{array}{rcl}
\nabla^{\Sigma}_{e_1} \sigma 
&=&    \partial_{e_1} \sigma + \frac{1}{4} \cdot 
\left\{
  \frac{d^2 T}{2} \gamma(E_2) \gamma(E_3) 
- \frac{d^2 T}{2} \gamma(E_3) \gamma(E_2)
\right\} \cdot \sigma \\
&=&    \partial_{e_1} \sigma + \frac{d^2 T}{4} 
       \left(\begin{array}{cc}  
       0 & i \\ i & 0 
       \end{array}\right) \cdot\sigma, \\
\nabla^{\Sigma}_{e_2} \sigma 
&=&    \partial_{e_2} \sigma + \frac{d^2 T}{4} 
       \left(\begin{array}{cc}  
       0 & -1 \\ 1 & 0 
       \end{array}\right) \cdot\sigma, \\ 
       \nabla^{\Sigma}_{e_3} \sigma 
&=&    \partial_{e_3} \sigma + \frac{d^2 T}{4} 
       \left(\begin{array}{cc}  
       -i & 0 \\ 0 & i 
       \end{array}\right) \cdot\sigma.     
\end{array}
$$
For the Dirac operator we obtain
$$
D=\sum\limits^3_{i=1} \gamma(e_i)\nabla^{\Sigma}_{e_i} = 
\sum\limits^3_{i=1} \gamma(e_i)\partial_{e_i} - \frac{d^2 T}{4}.
$$
Spin structures of $M(r,d,T)$ are given by homomorphisms 
$\epsilon:\Gamma_r \rightarrow \ZZ/2\ZZ = \{-1, 1\}$. 
Since $\ZZ/2\ZZ$ is abelian $\epsilon$ must factor through 
$\Gamma_r / [\Gamma_r, \Gamma_r]$. 
The commutator subgroup is 
$[\Gamma_r, \Gamma_r]=\{g(0,0,r\, z)|z\in\ZZ\}$. 
Hence 
$\Gamma_r / [\Gamma_r, \Gamma_r] \cong 
 r \, \ZZ\oplus\ZZ\oplus(\ZZ / r \, \ZZ)$. 
The homomorphism $\epsilon$ is determined by its image on generators of 
$\Gamma_r / [\Gamma_r, \Gamma_r]$, i.e.~by 
$$
\begin{array}{rcl}
\delta_1 &:=& \epsilon(g(r,0,0)),\\
\delta_2 &:=& \epsilon(g(0,1,0)),\\
\delta_3 &:=& \epsilon(g(0,0,1)).
\end{array}
$$
Put $\delta:=(\delta_1, \delta_2, \delta_3)$. 
The first two components $\delta_1$ and $\delta_2$ can take the values 
$1$ and $-1$. 
The third component $\delta_3$ can also take both values if $r$ is even, 
but if $r$ is odd $\delta_3$ is necessarily $\delta_3 = 1$. 
Summarizing, we see that spin structures on $M(r, d, T)$ are given 
by triples $\delta = (\delta_1, \delta_2, \delta_3) \in (\ZZ / 2 \ZZ)^3$ 
where $\delta_3$ must be $+1$ in case $r$ is odd. 
The corresponding homomorphism is then
$$
\epsilon (g(rx, y, z)) = \delta_1^x \cdot \delta_2^y \cdot \delta_3^z 
\in \ZZ / 2 \ZZ
$$
and spinor fields are maps $\sigma:G\rightarrow\Sigma_3=\CC^2$ such that 
\begin{eqnarray}
\sigma(g_0 g) = \epsilon(g_0)\cdot\sigma(g), \, 
g_0 \in \Gamma_r, \, g\in G. \label{heisenspinor}
\end{eqnarray}
\indent
The Heisenberg group $G$ acts on the Hilbert space $L^2(\Sigma_\ep M)$
of $L^2$-spinor fields on $M$ by the right regular representation $R$:
$$
\Big(R(g_0)\sigma\Big) (g) = \sigma (g\, g_0).
$$
At the same time the Clifford algebra $\Cl(3)$ acts on this Hilbert space 
by pointwise Clifford multiplication $\gamma$. 
Note that this action is compatible with condition (\ref{heisenspinor}). 
Moreover, the actions of $G$ and of $\Cl(3)$ commute. 
Modules with commuting $G$- and $\Cl(3)$-operations
will be called {\em $G$-$\Cl(3)$-bimodules}. 

Differentiation of a spinor field is given by the derived right regular 
action: 
$$
\partial_X \sigma = \stelle{\frac{d}{dt}}{t=0} R(\exp tX)\sigma =
R_{*} (X) \sigma, \quad X\in \mbox{{\goth g}}. 
$$

We will decompose the Hilbert space $L^2(\Sigma_\ep M)$ into closed 
subspaces invariant under the actions $R$ and $\gamma$. 
These subspaces are then also left invariant by the Dirac operator. 
Expressed in terms of $R$ and $\gamma$ the Dirac operator is given by 
\begin{eqnarray}
D=\sum\limits^3_{i=1} R_{*} (e_i) \otimes \gamma (E_i) - \frac{d^2 T}{4}.
\label{heisendirac}
\end{eqnarray}
\indent
Let the spin structure corresponding to 
$\delta = (\delta_1, \delta_2, \delta_3) \in (\ZZ / 2 \ZZ)^3$ be fixed. 
Let us start with the case $\delta_3 = 1$. 
Then for fixed $g\in G$ and spinor field $\sigma : G\rightarrow\Sigma_3$ 
the map 
$$
\varphi_g : \RR\rightarrow\Sigma_3, \, \varphi_g (z) = 
\Big(R(g(0,0,z))\sigma\Big)(g),
$$
is $1$-periodic, $\varphi_g(z+1)=\varphi_g(z)$. 
We expand $\varphi_g$ into a Fourier series 
$$
\varphi_g(z)=\sum\limits_{\tau\in\ZZ} \varphi_{\tau} (g) \cdot 
e^{2\pi i\tau z},
$$
where the Fourier coefficients are given by 
$\varphi_{\tau} (g) = \int\limits^1_0 R(g(0,0,t)) \sigma(g) 
e^{-2\pi i\tau t} dt$. 
Putting $z=0$ we obtain 
$$
\sigma(g) = \varphi_g (0) = \sum\limits_{\tau\in\ZZ} \varphi_{\tau}(g). 
$$
\indent
We found a first decomposition of $L^2(\Sigma_\ep M)= 
\bigoplus\limits_{\tau\in\ZZ} H_{\tau} = H_0 \oplus 
\bigoplus\limits_{\tau\in\ZZ -\{0\}} H_{\tau}$. 
On $H_{\tau}$ central elements $g(0,0,z)$ act via $R$ by multiplication 
by $e^{2\pi i\tau z}$. Hence $H_{\tau}$ decomposes into copies of the 
irreducible representations $\rho_{\tau}$. 
The subspace $H_0$ is the part on which the center of $G$ acts trivially. 
Thus it decomposes under $R$ into 1-dimensional representations of the 
form $\rho_{\alpha, \beta}$. 
We now determine this decomposition.

For $\sigma\in H_0$ the map $\RR^2 \mapsto\Sigma_3$, 
$(x,y)\rightarrow\sigma(g(x,y,0))$, is periodic for the lattice 
$2(r\ZZ\oplus\ZZ)$. 
Hence we can expand it into the Fourier series
$$
\sigma(g(x,y,z)) = \sum\limits_{
\begin{array}{c}
\alpha\in\frac{1}{2r}\ZZ \\
\beta\in\frac{1}{2}\ZZ
\end{array}
}
a_{\alpha, \beta} \cdot e^{2\pi i(\alpha x + \beta y)}.
$$
The equivariance property (\ref{heisenspinor}) imposes further
restrictions on $\alpha$ and $\beta$, namely 
$$
\sigma(g(x,y,z)) = \sum\limits_{
\begin{array}{c}
\alpha\in\frac{1}{2r}\ZZ, e^{2\pi i r \alpha} = \delta_1 \\
\beta\in\frac{1}{2}\ZZ, e^{2\pi i \beta} = \delta_2
\end{array}
}
a_{\alpha, \beta} \cdot e^{2\pi i(\alpha x + \beta y)}.
$$
This yields the decomposition 
$$
H_0 = \bigoplus\limits_{
\begin{array}{c}
\alpha\in\frac{1}{2r}\ZZ, e^{2\pi i r \alpha} = \delta_1 \\
\beta\in\frac{1}{2}\ZZ, e^{2\pi i \beta} = \delta_2
\end{array}
}
H^{\alpha, \beta}
$$
where each summand $H^{\alpha, \beta}$ is 
isomorphic to $\CC\otimes\Sigma_3 = \Sigma_3$ with $G$ acting via $R$ by 
$\rho_{\alpha, \beta}$ on the first factor 
(see \ref{rhoab} for the definition) and $\Cl(3)$ acting 
via $\gamma$ on the second. 

How does the Dirac operator $D$ (\ref{heisendirac}) act on $H^{\alpha, \beta}$? 
We compute 
$$
\begin{array}{rcl}
D &=& - d\, R_{*}(X)\cdot\gamma(E_1) - d \, R_{*}(Y)\cdot\gamma(E_2) \\
\\
  &&  +T^{-1} R_{*}(Z)\cdot\gamma(E_3) - \frac{d^2 T}{4} \\
\\
  &=& -d\cdot 2\pi i \alpha \cdot 
  \left(\begin{array}{cc}
  0 & i\\
  i & 0
  \end{array}
  \right)
  -d\cdot 2\pi i \beta \cdot 
  \left(\begin{array}{cc}
  0 & -1\\
  1 & 0
  \end{array} 
  \right) \\
\\
  && +T^{-1} \cdot 0 \cdot 
  \left(\begin{array}{cc}
  i & 0\\
  0 & -i
  \end{array}
  \right) 
  - \frac{d^2 T}{4} \cdot 
  \left(\begin{array}{cc}
  1 & 0\\
  0 & 1
  \end{array}
  \right) \\
\\
  &=& 
  \left(\begin{array}{cc}
  -\frac{d^2 T}{4} & 2\pi \mbox{\it id} (\beta - i \alpha) \\
  -2\pi \mbox{\it id} (\beta + i \alpha) & -\frac{d^2 T}{4}
  \end{array}
  \right).
\end{array}
$$
The eigenvalues are $-\frac{d^2 T}{4} \pm 2\pi d\sqrt{\alpha^2 + \beta^2}$. 

Next we decompose $H_{\tau}$, $\tau\neq 0$, under the actions of $G$ via $R$ and of 
$\Cl(3)$ via $\gamma$. 
Let $\sigma\in H_{\tau}$. 
Then $\sigma$ is determined by $f_{\sigma} : \RR^2 \rightarrow\CC$, 
for $f_{\sigma} (x,y)=\sigma(g(x,y,0))$, because 
$$
\sigma(g(x,y,z)) = e^{2\pi i \tau z} \cdot f_{\sigma} (x,y).
$$
Fixing $x\in\RR$ we expand $y\rightarrow f_{\sigma} (x,y)$ into a Fourier 
series and using (\ref{heisenspinor}) we obtain 
$$
f_{\sigma} (x,y) = \sum\limits_{\beta\in {1\over 2} \ZZ, 
e^{2\pi i \beta}=\de_2} 
b^{\sigma}_{\beta} (x) \cdot e^{2\pi i \beta y}.
$$
Using (\ref{heisenspinor}) once more we get the periodicity property of 
$b^{\sigma}_{\beta}$, 
$$
b^{\sigma}_{\beta}(x+\xi) = 
{\de_1}^{\xi/r} \cdot b^{\sigma}_{\beta+\tau\xi} (x), \, 
x\in\RR, \xi\in r \cdot\ZZ.
$$
Hence we have $|\tau| \cdot r$ independent functions 
$b^{\sigma}_{1 + \ep} ,\ldots, 
b^{\sigma}_{|\tau| \cdot r + \ep}$, with $\ep=1/2$ for $\de_2=-1$ and 
$\ep=0$ for $\de_2=1$.
The map 
\begin{eqnarray*}
H_{\tau} & \rightarrow & 
\bigoplus\limits^{|\tau|r}_{j=1} L^2(\RR, \Sigma_3) = 
\bigoplus\limits^{|\tau|r}_{j=1} L^2(\RR, \CC) \otimes \Sigma_3\\
\sigma & \mapsto & \left(b^{\sigma}_{1+\ep} ,\ldots, b^{\sigma}_{|\tau| r + \ep}\right)
\end{eqnarray*}
is an isomorphism. 
This decomposition is left invariant by the actions $R$ and $\gamma$. 
Each summand carries a $G$-$\Cl(3)$-bimodule structure given by the 
action of $R$ on the first factor and  
$\ga$ on the second. This $G$-$\Cl(3)$-bimodule structure is equivalent to
the $G$-$\Cl(3)$-bimodule structure given by $\rho_{\tau}$ and $\ga$ 
(for the definition of $\rho_\tau$ see \ref{rhot}).

We compute the eigenvalues of the Dirac operator on 
$L^2(\RR, \CC) \otimes \Sigma_3$. 
$$
\begin{array}{rcl}
D &=& -d(\rho_{\tau})_{*} (X) \cdot \gamma(E_1) -d(\rho_{\tau})_{*} (Y) 
\cdot\gamma(E_2) \\
\\ 
&& + T^{-1}(\rho_{\tau})_{*} (Z) \cdot \gamma(E_3) - \frac{d^2 T}{4} \\
\\
&=& -d \cdot \left(\begin{array}{cc}
  0 & i\\
  i & 0
  \end{array}
  \right) 
  \cdot \left(-\frac{d}{dt}\right) \, 
  -d \cdot \left(\begin{array}{cc}
  0 & -1\\
  1 & 0
  \end{array}
  \right) 
  (-2\pi i\tau t) \\
\\
&&  + T^{-1} \cdot \left(\begin{array}{cc}
  i & 0\\
  0 & -i
  \end{array}
  \right) 
  \cdot 2\pi i \tau \, -\frac{d^2 T}{4} \cdot 
  \left(\begin{array}{cc}
  1 & 0\\
  0 & 1
  \end{array}
  \right) \\
\\
&=& \left(\begin{array}{cc}
  -2\pi T^{-1} \tau \, -\frac{d^2 T}{4} 
  & -2\pi id \tau t + id \frac{d}{dt}\\
  2\pi id \tau t + id \frac{d}{dt} 
  & 2\pi T^{-1} \tau \, -\frac{d^2 T}{4}
  \end{array}
  \right).
\end{array}
$$
We use the $L^2$-basis $u_0, u_1, \ldots$ for $L^2(\RR, \CC)$ where 
$u_k(t) = h_k \left(\sqrt{2\pi|\tau|t} \, \right)$, see Section 2.
In case $\tau>0$ we apply $D$ to $u_k \choose 0$ and obtain
$$
\begin{array}{rcl}
D \mbox{$u_k \choose 0$} 
&=&
\left(
\begin{array}{l}
\left(-2\pi T^{-1}\tau \, -\frac{d^2 T}{4}\right) u_k
\\
\mbox{\it id} (2\pi\tau\,t\,u_k + u'_k)
\end{array}
\right) \\
\\
&\begin{array}{c}
(\ref{hermite1})\\=
\end{array}
&
\left(
\begin{array}{l}
\left(-2\pi T^{-1}\tau \, -\frac{d^2 T}{4}\right) u_k
\\
\mbox{\it id} (4\pi\tau\,t\,u_k + \sqrt{2\pi\tau} u_{k+1})
\end{array}
\right) \\
\\
&\begin{array}{c}
(\ref{hermite2})\\=
\end{array}
&
\left(
\begin{array}{l}
\left(-2\pi T^{-1}\tau \, -\frac{d^2 T}{4}\right) u_k
\\
- \mbox{\it id} \sqrt{2\pi\tau} \cdot 2k \, u_{k-1}
\end{array}
\right).
\end{array}
$$
Similarly, $D \mbox{$0 \choose u_{k-1}$} = 
\left(
\begin{array}{c}
\mbox{\it id} \sqrt{2\pi\tau} u_k \\
\left(2\pi T^{-1}\tau \, -\frac{d^2 T}{4}\right) u_{k-1}
\end{array}
\right)
$. 
Hence $u_0 \choose 0$ is an eigenvector for $D$ with eigenvalue 
$-2\pi T^{-1}\tau \, -\frac{d^2 T}{4}$ and for $k\in\NN$ the vectors 
$u_k \choose 0$ and $0 \choose u_{k-1}$ span a $D$-invariant subspace. 
With respect to this basis $D$ is given by the $2\times 2$-matrix
$$
D \sim 
\left(
\begin{array}{ll}
-2\pi T^{-1}\tau \, -\frac{d^2 T}{4} & \mbox{\it id} \sqrt{2\pi\tau}\\
-2 k \, \mbox{\it id}\sqrt{2\pi\tau} & 2\pi T^{-1}\tau \, -\frac{d^2 T}{4}
\end{array}
\right).
$$
The eigenvalues are 
$-\frac{d^2 T}{4}\pm 2\sqrt{\pi^2 \tau^2 T^{-2}+k\pi d^2 \tau}$. 
Analogous reasoning for the case $\tau<0$ yields an eigenvalue 
$-2\pi T^{-1} |\tau| \, -\frac{d^2 T}{4}$ and for $k\in\NN$ the 
eigenvalues 
$-\frac{d^2 T}{4} \pm 2 \sqrt{\pi^2 \tau^2 T^{-2}+k\pi d^2 |\tau|}$. 

It remains to consider the case $\delta_3 = -1$. 
In this case we get the decomposition 
$$
L^2 (\Sigma_\ep M)=\bigoplus\limits_{\tau\in\ZZ+\frac{1}{2}} H_{\tau},
$$
the representations $\rho_{\alpha, \beta}$ do not occur. 
The computation of the eigenvalues of $D$ on $H_{\tau}$ is the same as 
above. 
We summarize
\begin{theorem}\label{heisenspec}
Spin structures on the 3-dimensional Heisenberg manifold 
$M(r,d,T)$, $d$, 
$T>0$, $r\in\NN$, correspond to triples 
$\delta=(\delta_1, \delta_2, \delta_3) \in (\ZZ/2\,\ZZ)^3$ where the case 
$\delta_3 = -1$ occurs if and only if $r$ is even. 
The eigenvalues of the Dirac operator are as follows 
\begin{itemize}
\item[{\rm\bf A.}] If $\delta_3 = +1$:
\begin{itemize}
\item[{\rm a)}]
$\lambda^{\pm}_{\alpha, \beta} := -\frac{d^2 T}{4} \pm 2 \pi \, d 
\sqrt{\alpha^2 + \beta^2}$ 
\newline 
with multiplicity $1$ for each 
$\alpha\in\frac{1}{2r}\cdot\ZZ$ with 
$e^{2\pi\, i\, r\, \alpha} = \delta_1$ and each 
$\beta\in\frac{1}{2}\cdot\ZZ$ with $e^{2 \pi \, i \, \beta} = \delta_2$. 
\item[{\rm b)}]
$\lambda_{\tau} := -\frac{d^2 T}{4} \, - 2\pi\,\tau\,T^{-1}$ 
\newline
with multiplicity $2\tau r$ for each $\tau\in\NN$. 
\item[{\rm c)}]
$\lambda^{\pm}_{\tau, k} := -\frac{d^2 T}{4} \pm 2 
\sqrt{\pi^2 \tau^2 T^{-2} + k\, \pi d^2 \tau}$ 
\newline
with multiplicity $2\tau r$ for all $\tau$, $k\in\NN$. 
\end{itemize}
\item[{\rm\bf B.}] If $\delta_3 = -1$:
\begin{itemize}
\item[{\rm a)}]
$\lambda_{\tau} := -\frac{d^2 T}{4} \, - 2\pi\,\tau\,T^{-1}$ 
\newline 
with multiplicity $2\tau r$ for each $\tau\in\NN_0 +\frac{1}{2}$. 
\item[{\rm b)}]
$\lambda^{\pm}_{\tau, k} := -\frac{d^2 T}{4} \pm 2 
\sqrt{\pi^2 \tau^2 T^{-2} + k\, \pi d^2 \tau}$ 
\newline 
with multiplicity $2\tau r$ for each $\tau\in\NN_0 +\frac{1}{2}$ 
and each $k\in\NN$. 
\end{itemize}
\end{itemize}
\end{theorem}
Note that the Dirac spectrum for the spin structure 
$\delta=(\delta_1, \delta_2, -1)$ is independent of $\delta_1$ and of 
$\delta_2$. 

The eigenvalue $\lambda^+_{\tau, k}$ is zero for suitable $T>0$. 
Its multiplicity is $2\tau r$ hence unbounded if we let 
$\tau\rightarrow\infty$. 
Thus we have 

\begin{corollary}
For each $N\in\NN$ the 3-dimensional Heisenberg manifold 
$M = \Gamma_r\backslash G$ with spin structure corresponding to $\delta$ 
possesses a left invariant metric such that the dimension of the space of 
harmonic spinors is at least $N$.
\end{corollary}

A similar observation for Berger metrics on $S^3$ is due to Hitchin 
\cite{hi}. 
It is known that every closed 3-manifold with a fixed spin structure 
admits a Riemannian metric such that the corresponding Dirac operator has 
nontrivial kernel \cite{b2}. 
It is conjectured that it can be made arbitrarily large. 

The manifold $M(r, d, T)$ is a circle bundle over the torus 
$T^2 = \RR^2 / r \, \ZZ\oplus\ZZ$. 
If one shrinks the fiber length $(T\rightarrow 0)$ then $M(r, d, T)$ is 
said to {\it collapse} to $T^2$. 
It is interesting to examine the behavior of the eigenvalues under this 
collapse. 
We see that eigenvalues of the type $\lambda_{\tau}$ or 
$\lambda^{\pm}_{\tau, k}$ tend to $\pm \infty$. 
But the eigenvalues $\lambda^{\pm}_{\alpha, \beta}$ converge to 
$\pm 2\pi\,d \sqrt{\alpha^2 + \beta^2}$ which are precisely the Dirac 
eigenvalues of the torus $T^2$. 
\begin{corollary}\label{heisencollaps}
If the Heisenberg manifold $M(r, d, T)$ collapses to the 2-torus, i.e.~if 
$T\rightarrow 0$, then the Dirac eigenvalues behave as follows:

\begin{itemize}
\item[{\rm\bf A.}] If $\delta_3 = +1$:
\newline
There are eigenvalues converging to those of the 2-torus.
All other eigenvalues tend to $\pm\infty$. 
\item[{\rm\bf B.}] If $\delta_3 = -1$:
\newline
All eigenvalues tend to $\pm\infty$. 
\end{itemize}
\end{corollary}
This phenomenon will be studied in larger generality in the next section. 

To conclude this section we mention that a computation of the Dirac 
eigenvalues on Heisenberg manifolds can also be carried out in higher 
dimensions. 
It just becomes notationally more complicated, see \cite{b1} for more 
details. 
In all dimensions $\ge 5$ it is possible to find Heisenberg manifolds 
with nonisomorphic fundamental groups having the same Dirac spectrum.
See Section 5 for a general discussion.


\section{$S^1$-Actions and Collapse}

In this section we will study the asymptotic behavior of the 
spectrum of the Dirac operator on a larger class of collapsing manifolds.
This will generalize the observations we made in the last section 
for 3-dimensional Heisenberg manifolds.
Our discussion in this section is similar in spirit to the one
in \cite{bebebou} for the Laplace operator on functions.

We suppose that $S^1$ acts freely and isometrically on the 
compact connected  Riemannian spin manifold 
$(M,\ti g)$ of dimension $n+1$.
We view $M^{n+1}$ as the total space of an $S^1$-principal 
bundle over the base space $N^n=M^{n+1}/S^1$.

This base space shall carry the unique metric $g$ such that the projection
  $$ \pi : (M,\ti g) {\longrightarrow} (N,g)$$
is a Riemannian submersion.

This $S^1$-principal bundle has a unique connection 1-form
$i\om: TM \to i \RR$
such that $\ker \res{\om}{m}$ is orthogonal to the fibers for all $m \in M$
with respect to $\ti g$.

The $S^1$-action induces a Killing vector field $K$. 
To keep the discussion simple we will assume that the length $|K|$ is
constant on $M$.
This is equivalent to saying that the fibers of $\pi$ are totally
geodesic.
One can relax this assumption, compare Remark 4.2.

By rescaling the metric $\ti g$ on $M$ along the fibers while keeping it
the same on $\ker \om$ we obtain a 1-parameter family of metrics
$\ti g_\el$ on $M$ for which $\pi : M \to N$ is a Riemannian submersion,
$\om(K)\equiv 1$, and $2 \pi \el$ is the length of the fibers.
The length of $K$ with respect to $\ti g_\el$ is $\el$.

As the metric on $M$ is completely characterized by the 
connection 1-form $i\om$, the fiber length  $2 \pi \el $ 
and the metric $g$ on $N$, we can express the Dirac operator 
in terms of $\om$, $\el$, and $g$.
This will allow us to analyze the behavior of the spectrum 
for collapsing $S^1$-fibers $(\el \to 0)$. 
It will turn out that the spin structure on $M$
determines whether there are converging eigenvalues or not.

The $S^1$-action on $M$ induces an $S^1$-action on $\pso(M)$.
A spin structure $\ti\ph:\psp(M) \to \pso(M)$ will
be called {\it projectable} 
if this $S^1$-action on $\pso(M)$ lifts to $\psp(M)$. 
Otherwise it will be called {\it nonprojectable}.

Any projectable spin structure induces a spin structure on $N$ as follows:
Let $\ti\ph:\psp(M) \to \pso(M)$ denote the projectable spin structure on $M$ 
and let $\PsonM$ consist of those frames over $M$ having $K/\el $ as
the first vector.
We can identify $\pso(N)$ with $\PsonM/S^1$. 
Now $\ti\ph^{-1}(\PsonM)/S^1$ is a 
$\Spinn$-bundle over $N$ and $\ti\ph$ induces a corresponding spin structure
on $N$.

On the other hand, any spin structure on $N$ canonically induces a 
projectable spin structure on $M$ via pull-back:
Let $\ph:\psp(N)\to \pso(N)$ be a spin structure of $N$. 
Then $\pi^*\ph:\pi^*\psp(N)\to \pi^* \pso(N)=:\PsonM$ is a 
$\Theta_n:\Spin(n)\to \SO(n)$ equivariant map.
Enlarging the structure group to $\Spin(n+1)$ by
  $$\ti\ph:=\pi^*\ph\times_{\Theta_n}\Theta_{n+1}:
     \pi^*\psp(N)\times_{\Spin(n)}\Spin(n+1)\to \PsonM\times_{\SO(n)}\SO(n+1)$$
yields a spin structure on $M$.

Projectable spin structures and projectable spinors have already been studied
by Moroianu \cite{mo}.

Throughout this section we use the following convention: 
When we say that an operator has the eigenvalues $\mu_1,\mu_2,\ldots$
we suppose that each eigenvalue be repeated according to its multiplicity.  

\begin{theorem}\label{collapstheo}
Let $M$ be a closed Riemannian spin manifold.
Let $S^1$ act isometrically on $M$.
We assume that the orbits have constant length.
Let $N = M/S^1$ carry the induced Riemannian metric and let $\ti g_\el$
be the metric on $M$ described above.
Let $E \to N$ be a Hermitian vector bundle with a metric connection $\na^E$. 

We suppose that the spin structure on $M$ is projectable and that 
$N$ carries the induced
spin structure.
Let $\mu_1,\mu_2,\ldots$ be 
the eigenvalues of the twisted Dirac operator $D^E$ on $\Si N \otimes E \to N$.

Then we can number the eigenvalues 
$(\la_{j,k}(\el ))_{j \in \NN, k\in \ZZ}$ of the twisted Dirac
operator ${\ti{D}}^\el$ for $\ti g_\el$ on $\Si M \otimes \pi^*E \to M$
such that they depend continuously on $\el $ and such that for $\el \to 0$:
\begin{enumerate}[(1)]
\item For any  $j\in\NN$ and $k\in\ZZ$ 
         $$\el \cdot \la_{j,k}(\el) \to k. $$
      In particular,  
      $\la_{j,k}(\el ) \to \pm \infty$ if $k \not= 0$.
\item If $n=\dim N$ is even, then 
      $$\la_{j,0}(\el ) \to \mu_j.$$
      On the other hand if $n=\dim N$ is odd, then
      \begin{eqnarray*}
         \la_{2j-1,0}(\el ) & \to & \phantom{-} \mu_j\cr
         \la_{2j,0}(\el )   & \to &           - \mu_j
      \end{eqnarray*}
      In both cases, the convergence of the 
      eigenvalues $\la_{j,0}(\el )$ is uniform in $j$.
\end{enumerate} 
\end{theorem}

\begin{remark}
\em
The theorem also holds for a 
family of nonconstant length functions $\el :N\to {\RR}^+$  
provided $\el $  and $(1/\el )\grad \el $ converge uniformly to $0$.
The proof is analogous to the one below except that some formulas 
contain an additional term depending on $\grad \el $.
\end{remark}

In the proof of the theorem we will assume for simplicity
that $E$ is a trivial line bundle. 
For nontrivial $E$ every spinor bundle on $N$ or $M$ 
simply has to be twisted by $E$ or $\pi^*E$ resp. 
 
In order to prove Theorem \ref{collapstheo} we write the Dirac operator 
$\ti{D}^\el $ as a sum of a vertical Dirac operator, a horizontal 
Dirac operator, and a zero order term. 

For the definition of the vertical Dirac operator we need the Lie derivative
of a spinor along the $S^1$-fibers.
The action of $S^1$ on $\psp{M}$ induces an  
action of $S^1$ on $\Si M=\psp(M^{n+1})\times_{\Spin(n+1)}\Si_{n+1}$ 
which we denote by $\ka$. 
A spinor with base point $m$ will be 
mapped by $\ka(e^{it})$ to a spinor with base point $m\cdot e^{it}$.
We define the {\em Lie derivative} of a smooth spinor $\Psi$ 
in the direction of the Killing field
$K$ by
  $$\LL_K(\Psi)(m)= \stelle{d \over dt}{t=0}
    \ka(e^{-it})(\Psi(m\cdot e^{it})).$$   
See \cite{bougau} for a general discussion of Lie derivatives for spinors.
Since $\LL_K$ is the differential of a representation of the Lie group $S^1$ 
on $L^2(\Si M)$, we get the decomposition
  $$L^2(\Si M)=\bigoplus_{k \in \ZZ} V_k$$
into the eigenspaces $V_k$ of $\LL_K$ for the eigenvalue $ik$, $k\in \ZZ$.
The $S^1$-action commutes with the Dirac operator on $M$, 
hence this decomposition is preserved by the Dirac operator.  

We now want to calculate the difference between the covariant derivative and 
the Lie derivative in the direction of $K$. 
We always use the convention that any $r$-form 
$\al$ acts on a spinor $\Psi$ by
  $$\ga(\al)\Psi:=\sum_{i_1<\dots<i_r}\al(e_{i_1},\dots,e_{i_r})\,
    \ga(e_{i_1})\cdots \ga(e_{i_r}) \Psi$$
where the $e_i$ form an orthonormal basis of the tangent space.

\begin{lemma}\label{vertinabla}
If $\ti\Psi$ is a smooth section of $\Si M$, then 
  $$ \na_K \ti\Psi -  \LL_K \ti\Psi= 
    {\el^2 \over 4}\, \ga(d\om) \ti\Psi$$ 
\end{lemma}

\proof{}
First we fix some notations. 
In the rest of this section, we denote
the horizontal lift of $X \in TN$ by 
$\ti{X}\in TM$.
Let $(f_1,\dots, f_n)$ be a local orthonormal frame of $N$ on the open subset 
$U$ and let 
$A$ be a lift to the spin structure. 
Then $(e_0:=K/\el ,e_1:=\ti{f}_1\circ \pi,\dots,e_n:=\ti{f}_n\circ\pi)$ 
is a local orthonormal frame on $\pi^{-1}(U)$ with lift $\pi^*A$.
Using the Koszul formula we can express the 
Christoffel symbols $\GaM{ij}{k}$ for 
$(e_0,\dots,e_n)$ in terms of the Christoffel symbols $\GaN{ij}k$ 
of $(f_1,\dots,f_n)$,
the $S^1$-bundle curvature form $d\om$ and the length $\el $.
For $i,j,k \in \{1,\ldots,n\}$ we get 
\begin{equation}\label{gammamformel}
  \GaM{ij}k=\GaN{ij}k \qquad 
  -\GaM{ij}0=\GaM{i0}j=\GaM{0i}j= {\el \over 2}\, d\om (e_i,e_j) 
\end{equation}
  $$\GaM{i0}0=\GaM{0i}0=\GaM{00}i=\GaM{00}0=0$$
With respect to these local frames on $M$ and $N$ 
we apply (\ref{nablaspinformel}) for $M$ and we obtain
\begin{eqnarray*}
 {1\over \el } \na_K \ti\Psi & = 
    & {1\over \el } \pa_K\ti\Psi + {1\over 4}
      \sum_{j,k=0}^n{\GaM{0j}k}\ga(e_j)\ga(e_k) \ti\Psi\\
  &=& {1\over \el } \pa_K\ti\Psi +{\el  \over 8} \sum_{j,k=1}^n 
      d\om(e_j,e_k)\ga(e_j)\ga(e_k) \ti\Psi\\  
  &=& {1\over \el } \pa_K\ti\Psi +{\el  \over 4}\,\ga(d\om) \ti\Psi.    
\end{eqnarray*}
On the other hand the choice of the frames implies that 
  $$\LL_K\ti\Psi=\pa_K\ti\Psi,$$
thus the lemma follows.
\qed 
 
Now we associate to the $S^1$-bundle $M\to N$ the complex line bundle 
$L:=M\times_{S^1}\CC$ with the natural connection given by $i \om$.

\begin{lemma}\label{horinabla}
Let $n$ be even.
There is a homothety of Hilbert spaces 
  $$Q_k:L^2(\Si N \otimes L^{-k}) \to  V_k$$
such that the horizontal covariant derivative is given by
  $$\na_{\ti X}{Q_k(\Psi)} =Q_k(\na_{X}\Psi) 
    + {\el \over 4}\, \ga(K/\el )\ga(\ti V_X){Q_k(\Psi)}$$
where $V_X$ is the vector field on $N$ satisfying
$d\om(\ti X,\cdot)=\langle \ti V_X, \cdot \rangle$.
Clifford multiplication is preserved, i.e.\   
   $$Q_k(\ga(X)\Psi)=\ga(\ti X) {Q_k(\Psi)}.$$
\end{lemma} 

\proof{}
Since $n$ is even, we have $\Si_n=\Si_{n+1}$.
We define a vector bundle map
\begin{eqnarray*}
  \Pi_k: \Si M &\to & \Si N \otimes L^{-k}\\
  (m,[\pi^*A,\si]) & \mapsto & (\pi(m), [A,\si]\otimes [m,1]^{-k}),
\end{eqnarray*}
where $[\pi^*A,\si]$ (resp.\ $[A,\si]$, $[m,1]$) denotes the equivalence
class of $(\pi^*A,\si)$ (resp.\ $(A,\si)$, $(m,1)$) in 
$\pi^*\psp(N)\times_{\Spin(n)}\Si_n=\Si M$ 
(resp.\ $\psp(N)\times_{\Spin(n)}\Si_n=\Si N$, resp.\ $M\times_{S^1}\CC=L$). 
This bundle map is fiberwise a vector space isomorphism preserving Clifford 
multiplication. 
Therefore for any spinor $\Psi:N \to \Si N \otimes L^{-k}$ there is a unique
${Q_k(\Psi)}$ such that

\vspace{-1cm}
$$
\renewcommand{\arraystretch}{5.5}
\hfill
\begin{array}{r@{\hskip 2cm}c}
\rnode{a}{M} & \rnode{b}{\Si M} \\
\rnode{c}{N}&\rnode{d}{\Si N \otimes L^{-k}}   \\
\end{array}
\psset{nodesep=3pt}
\ncline{->}{a}{b}\Aput{{Q_k(\Psi)}}
\ncline{->}{a}{c}\Aput{\pi}
\ncline{->}{b}{d}\Aput{\Pi_k}
\ncline{->}{c}{d}\Aput{\Psi}
\hfill
$$
commutes.
Hence $Q_k$ is a well defined homomorphism of Hilbert spaces and injectivity
of $Q_k$ follows from the surjectivity of $\pi$.

For any section $\ti{\Psi}$ of $\Si M$ we get the equivalences
\begin{eqnarray*}
& \ti{\Psi} \in \image(Q_k)\\
\gdw& \Pi_k \circ \ti\Psi (m)= \Pi_k \circ \ti\Psi(me^{it})& 
\forall m\in M, e^{it} \in S^1 \\
\gdw& \ti\Psi(m)=e^{-ikt} \ka(e^{-it}) \ti\Psi(m e^{it})& 
\forall m\in M, e^{it} \in S^1 \\
\gdw& ik\, \ti\Psi= \LL_K \ti\Psi.
\end{eqnarray*} 
Thus the image of $Q_k$ is precisely $V_k$. 
The map $Q_k$ is an isometry up to the factor~$\el$.   
 
Finally it remains to show the formula for the covariant derivative in the 
horizontal directions. For this it suffices to prove it locally 
and we can restrict to the case $X=f_i$ with the notations as in the 
previous lemma. 
As before we use (\ref{nablaspinformel}) and (\ref{gammamformel}).
\begin{eqnarray*}
  \na_{{\ti f}_i} {Q_k(\Psi)} & = 
  & Q_k(\na_{f_i}\Psi) + {1 \over 4} \sum_{j=1}^n
    \GaM{i0}j\ga(K/\el )\ga(e_j){Q_k(\Psi)}\\
  && + {1\over 4} \sum_{j=1}^n \GaM{ij}0 \ga(e_j)\ga(K/\el ) {Q_k(\Psi)}\\
  &=& Q_k(\na_{f_i}\Psi) + {\el \over 4}\, 
     \ga(K/\el )\ga(\ti V_{f_i}){Q_k(\Psi)}
\end{eqnarray*}
\qed

\proof{of Theorem \ref{collapstheo} for $n$ even.}
We define the {\it horizontal Dirac operator} as the unique closed linear
operator $D_h:L^2(\Si M) \to L^2(\Si M)$ 
on each $V_k$ given by
  $$D_h:= Q_k\circ D \circ {Q_k}^{-1}$$
where $D$ is the (twisted) Dirac operator on $\Si N \otimes L^{-k}$.

We define the {\it vertical Dirac operator} 
  $$D_v:=\ga(K/\el)\,\LL_K$$ 
and the zero order term 
  $$Z:=-(1/4)\,\ga(K/\el )\,\ga(d\om).$$
 
Using Lemma \ref{vertinabla} and Lemma \ref{horinabla} we can express 
the Dirac operator $\ti D^\el $ as a sum
\begin{eqnarray*}
   \ti D^\el & = & \sum_{i=0}^n \ga(e_i)\na_{e_i} \\
          & = & {1\over \el } D_v + D_h + \el  Z.\\
\end{eqnarray*}
Since $D_h$, $\ga(K/\el)$ and $Z$ commute with the $S^1$-action they also 
commute with $\LL_K$. 
On the other hand, we have 
$\ga(K/\el)Q_k(\Psi)=Q_k(c \ga(\dvol_n)\Psi)$ with $c \in \{1,i,-1,-i\}$ 
depending on $n$ and the representation of $\Cl_{n+1}$.
Since $\ga(\dvol_n)$ anticommutes with any twisted Dirac operator on
$N$, we know that $\ga(K/\el)$ anticommutes with $D_h$.
Let $\Psi$ be a common eigenspinor for $\LL_K$ and $D_h$ for the 
eigenvalues $ik$ and $\mu$ resp.

On $U:=\span\{\Psi, \ga(K/\el ) \Psi\}$
the operator $A_\el := (1/\el) D_v + D_h$ is represented by the matrix
  $${1\over \el}\pmatrix{0& -ik\cr ik & 0} + \pmatrix{\mu & 0 \cr 0 & -\mu}=
    \pmatrix{\mu & -ik/\el  \cr ik/\el  & -\mu}.$$
 
Thus for $k=0$ the restriction $\res{A_\el }{U}$ has eigenvalues
$\pm \mu$. For $k\neq 0$ the eigenvalues of $\res{A_\el }{U}$ are 
the square roots of $(k/\el )^2+\mu^2$.
Therefore the eigenvalues $(\la_{j,k}^0(\el ))_{j\in\NN,k \in \ZZ}$ of 
$A_\el $ can be numbered such that they are continuous in $\el $ and 
satisfy properties (1) and (2) of Theorem \ref{collapstheo}. 
The additional term $\el Z$ does not change this behavior because
$Z:L^2(\Si M) \to L^2(\Si M)$ is bounded.
\qed

\proof{of Theorem \ref{collapstheo} for $n$ odd.}
In the case ``$n$ odd'' the proof is similar but we have to do 
some modifications as the vector space equality $\Si_{n+1}=\Si_n$ 
is no longer valid.

For the standard basis $E_1,\dots,E_n \in \RR^n$ and $E_0,\dots,E_n
\in \RR^{n+1}$ we define the {\em volume element} by
  $$\nu_n:=i^{n+1\over 2} \ga(E_1)\cdots\ga(E_n) \in \Cl(n+1).$$
Its square is the identity and hence has the eigenvalues $\pm 1$. Let 
  $$\Si_{n+1}=\Si_n^{(+)} \oplus \Si_n^{(-)}$$
be the splitting into the eigenspaces.
These $\Si_n^{(\pm)}$ are the two irreducible representations of $\Cl(n)$.
The Clifford multiplication $\ga(E_0)$ anticommutes with $\nu_n$ and therefore
restricts to an isomorphism of $\Si_n^{(+)}\to \Si_n^{(-)}$. This isomorphism
changes the sign of $\ga(E_i)$ for $i>0$ as 
$\ga(E_0)\ga(E_i)=-\ga(E_i)\ga(E_0)$.
 
Denote the corresponding spinor bundles by $\Si^{(\pm)}N$. 
Lemma \ref{horinabla} is now valid with
  $$Q_k:L^2\Big((\Si^{(+)}N\oplus\Si^{(-)}N)\otimes L^{-k}\Big)\to V_k.$$
The same arguments as in the even case can be used to prove 
that the eigenvalues of $\res{\ti D^{\el}}{V_k}$ converge to $\pm \infty$
for $k\neq 0$ and to the eigenvalues of the Dirac operator on 
$\Si^{(+)}N\oplus\Si^{(-)}N$ for $k=0$.  
The Dirac operator on $\Si^{(+)}N$ has the same eigenvalues as the 
Dirac operator on $\Si^{(-)}N$ but with opposite signs.
This yields the statement of Theorem \ref{collapstheo}.
\qed

In the case of a nonprojectable spin structure we get a similar result.
Here the variable $k$ takes no longer integer values but 
  $$k \in \ZZ+{1\over 2}\, .$$ 

\begin{theorem}\label{collapsnoproj}
Let $M$ be a closed Riemannian spin manifold.
Let $S^1$ act isometrically on $M$.
We assume that the orbits have constant length.
Let $N = M/S^1$ carry the induced Riemannian metric and let $\ti g_\el$
be the metric on $M$ described above.
Let $E \to N$ be a Hermitian vector bundle with a metric connection $\na^E$. 

We suppose that the spin structure on $M$ is nonprojectable.

Then we can number the eigenvalues 
$(\la_{j,k}(\el ))_{j \in \NN, k\in \ZZ+(1/2)}$ of the twisted Dirac operator 
${\ti{D}}^\el$ for $\ti g_\el$ on $\Si M \otimes \pi^*E \to M$
such that they depend continuously on $\el $ and such that
for $\el  \to 0$
      $$\el \cdot \la_{j,k}(\el) \to k\quad \forall j=1,2,\dots.$$
In particular, $\la_{j,k}(\el ) \to \pm \infty$.
\end{theorem}

\proof{}
The proof is a variation of the proof of Theorem \ref{collapstheo}.
We will restrict to the case $n$ even and $E$ trivial.

Let $\ti\ph:\psp(M)\to \pso(M)$ be a nonprojectable spin structure on $M$. 
In this case $N$ may or may not be spin. 
But as before the principal 
$\SO(n+1)$-bundle $\pso(M)$ can be restricted to a principal $\SO(n)$-bundle 
$P_{\SO(n)}(M)$. 
Moreover $P:=\ti\phi^{-1}(P_{\SO(n)}(M))$ is a principal $\Spin(n)$-bundle 
over $M$. 

The action of $S^1\cong (\RR / 2\pi \ZZ)$ does not lift to $P$, 
but the double covering of $S^1$, i.e. $S^1\cong (\RR/ 4\pi \ZZ)$,
does act on $P$.
We define $\Spin^{\CC}(n)$ to be $\Spin(n)\times_{\ZZ_2}S^1$ 
where $-1\in\ZZ_2$ identifies $(-A,c)$ with $(A,-c)$. 
The Lie group $\Spin^{\CC}(n)$ naturally acts on $\Si_n$.
The actions of $\Spin(n)$ and $S^1$ on $P$ induce a free action of 
$\Spin^{\CC}(n)$ on $P$ and we can view $P$ as a principal 
$\Spin^{\CC}(n)$-bundle over $N$. Then we can form the bundle 
$P\times_{\Spin^{\CC}(n)} \Si_n$. 

If $N$ is spin, this bundle is just $\Si N \otimes L^{1\over 2}$. 
If $N$ is not spin, then neither $\Si N$ nor $L^{1\over 2}$ exist
but $\Si N \otimes L^{1\over 2}$ does exist.

Again we get a splitting 
  $$L^2(\Si M)=\bigoplus_{k\in \ZZ+{1\over 2}} V_k$$
into eigenspaces $V_k$ for $\LL_K$ to the eigenvalue $ik$. 
The rest of the proof of Theorem \ref{collapsnoproj} is the same as
the one for the case $k\neq 0$ in Theorem \ref{collapstheo}.
\qed

\noindent
{\bf Example}{}
We apply Theorems \ref{collapstheo} and \ref{collapsnoproj} to the 
Hopf fibration $S^{2m+1} \to \CC P^m$.
The Dirac eigenvalues for the metric $\ti g_\el$ on $S^{2m+1}$
have been computed by the second author in \cite[Theorem 3.1]{b2}:

{\it The Dirac operator of} $S^{2m+1}$ {\it with Berger
metric} $\ti g_\el$ {\it has the following eigenvalues}:
\begin{itemize}
\item[(i)]
$\lambda(\el) = \frac{1}{\el} \left( a+\frac{m+1}{2} \right) + \frac{\el m}{2}$, 
$a\in\NN_0$, {\it with multiplicity} $m+a \choose a$.
\item[(ii)]
$\lambda(\el) = (-1)^{m+1} \cdot 
\left\{
\frac{1}{\el} \left( a+\frac{m+1}{2} \right) + \frac{\el m}{2}
\right\}$, $a\in\NN_0$, {\it with multiplicity} $m+a \choose a$.
\smallskip
\item[(iii)]\ \newline\vskip-2.35cm
\begin{eqnarray*}
  \lambda_{1, 2}(\el) = (-1)^j \, \frac{\displaystyle \el}{\displaystyle 2} 
  & \pm & \Bigg\{
  \left[ \frac{\el}{2} (m - 2j - 1) + \frac{1}{\el} 
  \left( a_1 - a_2 + \frac{m-1}{2} - j\right) \right]^2 \\ 
  && \phantom{\Bigg\{} + 4 (m-j+a_1) (j+1+a_2)\Bigg\}^{1/2}
\end{eqnarray*}
\komment{
\item[(iii)]
$
\lambda_{1, 2}(\el) = (-1)^j \, \frac{\displaystyle \el}{\displaystyle 2} \pm 
\sqrt{ 
\begin{array}{l}
\left[ \frac{\el}{2} (m - 2j - 1) + \frac{1}{\el} 
\left( a_1 - a_2 + \frac{m-1}{2} - j\right) \right]^2 \\ 
 + 4 (m-j+a_1) (j+1+a_2)
\end{array}
}
$
}
\item[]
$a_1$, $a_2 \in \NN_0$, $j=0, \ldots, m-1$, {\it with multiplicity} 
$$
\frac{(a_1 + m) !\, (a_2 + m) !\, (a_1 + a_2 + m + 1)}{
a_1 !\, a_2 !\, (a_1 + m - j)\, (a_2 + j + 1)\, m !\, j !\, (m - j - 1)!}.
$$
\end{itemize}

\noindent
We observe for the eigenvalues of type (i) and (ii)
$$ 
\lim_{\el\to 0} \el\cdot\lambda(\el) = \pm\left( a+\frac{m+1}{2} \right).
$$ 
We see that the spin structure on $S^{2m+1}$ is projectable if and
only if $m$ is odd.
In fact, $\CC P^m$ is spin if and only if $m$ is odd.

To obtain a finite limit of $\lambda(\el)$ itself we must necessarily
have an eigenvalue of type (iii) and $a_1 - a_2 + \frac{m-1}{2} = j$.
In particular, $m$ must be odd again.
Plugging this expression for $j$ into the formula for the eigenvalue
of type (iii) and passing to the limit $\el \to 0$ yields the eigenvalues
of the complex projective space:

\begin{theorem}\label{cpnspektrum}
The Dirac operator on $\CC P^m$, $m$ odd, has the eigenvalues
$$
\pm 2
\sqrt{\left(a_1+\frac{m+1}{2}\right)\left(a_2+\frac{m+1}{2}\right)}
$$
where $a_1, a_2 \in \NN_0, |a_1-a_2| \le \frac{m-1}{2}$ with
multiplicity
$$
\frac{(a_1 + m) !\, (a_2 + m) !\, (a_1 + a_2 + m + 1)}{
a_1 !\, a_2 ! \left(a_1 + \frac{m+1}{2}\right) \left(a_2 + 
\frac{m+1}{2}\right) m ! \left(a_1-a_2+ \frac{m-1}{2}\right) !\, 
\left(a_2-a_1+\frac{m-1}{2}\right)!}.
$$
\end{theorem}

The formula for the eigenvalues coincides with the results in 
\cite{cfg1,cfg2,sese} after the right parameter substitutions.
Our expression however is simpler than the previous ones.
The formula for the multiplicities is new.


\section{Isospectral Deformations}

We now study the isospectrality question for the Dirac operator
on Riemannian nilmanifolds $\Gamma \backslash G$. 
The spectrum of the Dirac operator does not only depend on the Riemannian 
structure but also on the spin structure, i.e.\ on the homomorphism
$\ep:\Ga \to \{-1, 1\}$. 
For $\ep\equiv 1$ it will turn out in the examples under consideration
that the spectrum of 
the Dirac operator behaves much like the spectrum of the Laplacian on 1-forms.
For other $\ep$ the behavior can be different.

At first, let $L^2_{\ep}(\Ga \bs G)$ be the space of all complex-valued
locally square-integrable
functions on $G$, satisfying the $\ep$-equivariance condition
  $$f(g_0 g)=\ep(g_0) f(g),$$
for all $g_0 \in \Gamma$ and $g \in G$.
The nilpotent group
$G$ acts on this space via the right regular representation 
\begin{eqnarray*}
  R : G  &  \to &  \End(L^2_{\ep}(\Ga \bs G))\\
  g_0 & \mapsto &\big(R(g_0): f \mapsto f(\phantom{}\cdot\phantom{}
                                  g_0)\big).
\end{eqnarray*}
Let $\Ga'$ be another cocompact discrete subgroup of $G$,
$\ep':\Ga'\to\{-1, 1\}$ a homomorphism.
We say that $(\Ga,\ep)$ and $(\Ga',\ep')$ are {\it representation equivalent} 
if the right regular representations on $L^2_{\ep}(\Ga \bs G)$ and on
$L^2_{\ep'}(\Ga' \bs G)$
are isomorphic as $G$-representations.

As already mentioned at the end of Section 1, sections of the spinor 
bundle $\Si_\ep (\Ga \bs G)$ may be identified with $\Ga$-equivariant 
maps $G \to \Si_n$.
Elements of $\Si_n$ will also be interpreted as constant maps $G \to \Si_n$ 
and therefore as sections of $\Si G$.

Hence $f \otimes s \mapsto f\cdot s$ defines an isomorphism of Hilbert spaces 
$ L^2_{\ep}(\Ga \bs G) \otimes_{\CC} \Si_n \to L^2(\Si_\ep (\Ga \bs G))$.
In the rest of this section  we will identify via this isomorphism.
Clifford multiplication is performed on the second factor 
$$ 
\ga(X) ( f\otimes s)= f \otimes ( \ga(X) s)
$$  
for all $X \in \gg, f \in L^2_{\ep}(\Ga \bs G), s \in \Si_n$.

\begin{theorem}\label{repequi}
If $(\Ga ,\ep)$ and $(\Ga' ,\ep')$ 
are representation equivalent, then the corresponding Dirac operators 
have the same spectrum for any left-invariant Riemannian metric on $G$.
\end{theorem}

\noindent
{\bf Examples:}
\begin{itemize}
\item Milnor's 16-dimensional tori \cite{miln} 
satisfy the condition of the theorem for $\ep\equiv\ep'\equiv 1$.
\item The Gordon-Wilson deformations \cite{gw1} preserve the spectrum 
of the Dirac operator. 
In fact, if $\Psi_t$ is a family of 
almost inner automorphisms of $G$, $\Ga$ a cocompact discrete subroup
and $\ep:\Ga \to \{-1,1\}$ a homomorphism, then $(\Psi_t(\Ga),
\ep\circ{\Psi_t}^{-1})$ are pairwise representation equivalent.
This deformation preserves the length spectrum, the spectrum of the 
Laplacian on functions and on forms, and the Dirac spectrum for all
spin structures. 
\end{itemize}

Since the spectrum of the Dirac operator determines the volume, we get
\begin{corollary}\label{volcor}
If $(\Ga, \ep)$ and $(\Ga', \ep')$ 
are representation equivalent, then $\Ga\bs G$ and $\Ga'\bs G$ have 
the same volume for any left invariant volume form on $G$.
\end{corollary}

\proof{of Theorem \ref{repequi}}
Since $L^2_{\ep}(\Ga \bs G)$ and 
$L^2_{\ep'}(\Ga' \bs G)$ 
are isomorphic as representations of the Lie group $G$
they are isomorphic as representations of its Lie algebra $\gg$.
Let 
$$
\th: L^2_{\ep}(\Ga \bs G)\to L^2_{\ep'}(\Ga' \bs G)
$$
be such an isomorphism, i.e.\ $\th \circ R= R' \circ \th$ and 
$\th \circ R_* = {R'}_* \circ \th$.
In order to prove the theorem we show that $\th \otimes \id$ maps 
eigenspinors of the Dirac operator in $L^2(\Sigma_\ep(\Ga\bs G))$ to 
eigenspinors of the Dirac operator in $L^2(\Sigma_{\ep'}(\Ga' \bs G))$. 
It is important 
to keep in mind that these two Dirac operators are just restrictions of the 
Dirac operator on $\Si G$ to spinors equivariant under 
$\ep$ (resp.\ $\ep'$). 

Let $s=\sum_i f_i  s_i$ be an arbitrary eigenspinor in 
$L^2(\Si_\ep(\Ga \bs G))$  for the 
eigenvalue $\la$ with $f_i \in L^2_{\ep}(\Ga \bs G)$ and $s_i \in \Si_n$.
Let $E_1,\dots,E_n$ an orthonormal basis of $\gg$. 
Using the Leibniz rule
\begin{equation}\label{leibniz}
  \na_{E_j}(f_i s_i) = ( R_*(E_j)f_i)s_i + f_i\na_{E_j}s_i
\end{equation}
we get
\begin{eqnarray*}
\la (\th\otimes \id) (s) & = & (\th\otimes\id)(Ds)\\
&=&  \sum_{ij}(\th\otimes\id)\left( (\id \otimes \ga(E_j))
 \na_{E_j}(f_i s_i)\right)\\
& = & \sum_{ij} \left(\th (R_*(E_j)f_i)\ga(E_j)s_i+\th(f_i)\ga(E_j)
\na_{E_j}s_i\right)\\
&=& \sum_{ij} \left( {R'}_*(E_j)(\th(f_i)) \ga(E_j)s_i +\th(f_i)\ga(E_j)
\na_{E_j}s_i\right)\\
&=& \sum_{ij} (\id \otimes \ga(E_j)) \na_{E_j} (\th (f_i) s_i)\\
&=& D((\th \otimes \id) (s)).
\end{eqnarray*}
Thus $\th \otimes \id$ is an isomorphism of Hilbert spaces that preserves 
eigenspinors and eigenvalues and therefore the Dirac spectra coincide.
\qed

\komment{$M_i:=\Ga_i \bs G$ is an $S^1$-bundle over
 $\ol{M}_i:=\ol{\Ga}_i\bs G$.  
}

To get further applications of Theorem \ref{repequi} we need criteria 
telling us under which conditions $(\Ga ,\ep)$ and $(\Ga' ,\ep')$ 
are representation equivalent.
In the examples of Theorem \ref{repequi} we already mentioned that almost
inner automorphisms give us representation equivalent families.
This criterion is sufficient but not necessary. 
However, in the special case when $G$ is stricly nonsingular, 
we will get a necessary and sufficient criterion.
A nilpotent Lie group $G$ with center $Z(G)$
will be called {\it stricly nonsingular} if 
for every $z \in Z(G)$ and $x\in G\ohne Z(G)$ there exists 
an $a \in G$ such that 
  $$ xax^{-1}a^{-1} = z.$$
Note that the Heisenberg group is strictly nonsingular.

Let $\Ga$ be a cocompact discrete subgroup of $G$.
The quotient $G/Z(G)$ will be denoted by $\ol{G}$ and the image of $\Ga$
under the canonical projetion $\pi:G \to \ol G$ by $\ol{\Ga}$.
According to \cite[Proposition 2.17]{r} and \cite[Lemma 5.1.4]{cg} 
$\ol\Ga$ is a cocompact discrete subgroup of $\ol G$. 
If $\res\ep{\Ga\cap Z(G)}\equiv 1$ then there is an $\ol\ep$ such that
\vspace{-1cm}
$$
\renewcommand{\arraystretch}{5.5}
\begin{array}{r@{\hskip 2cm}l}
\rnode{a}{\Ga} & \rnode{b}{\{-1,1\}} \\
\rnode{c}{\ol\Ga}&  \\
\end{array}
\psset{nodesep=3pt}
\ncline{->}{a}{b}\Aput{\ep}
\ncline{->>}{a}{c}
\ncline{->}{c}{b}\Bput{\ol\ep}
$$
commutes. Such an $\ep$ will be called {\it projectable}. 
For $\dim Z(G)=1$ this definition agrees with our definition of 
projectable spin structures in Section 4.

\begin{theorem}\label{nonsingquot}
Let $G$ be a simply connected, strictly nonsingular nilpotent 
Lie group with left invariant metric. Let $\Ga_1$, $\Ga_2$ be 
cocompact, discrete subgroups of $G$ such that 
  $$\Ga_1\cap Z(G)=\Ga_2\cap Z(G)\hspace{.4cm}(=: \Ga_Z)$$
and let $\ep_i:\Ga_i \to \{-1,1\}$ be homomorphisms.
We define the volume quotient
$$v:={ \vol(\ol{\Ga}_1\bs\ol{G}) \over \vol(\ol{\Ga}_2\bs\ol{G})}
      ={\vol(\Ga_1\bs G)\over \vol(\Ga_2\bs G)}\in \QQ^+.$$
Then $(\Ga_1,\ep_1)$ and $(\Ga_2,\ep_2)$ are representation equivalent if
and only if one of the following conditions is satisfied:
\begin{enumerate}[(1)]
\item $\res{\ep_1}{\Ga_Z}\equiv\res{\ep_2}{\Ga_Z}\equiv 1$  
and $(\ol{\Ga}_1,\ol\ep_1)$ and $(\ol{\Ga}_2,\ol\ep_2)$ are 
representation equivalent for $\ol{G}$. 
\item $\res{\ep_1}{\Ga_Z}\equiv\res{\ep_2}{\Ga_Z}\not\equiv 1$ and $v=1$.
\end{enumerate}
\end{theorem}

\begin{remark}\label{multremark}\rm
If under the assumptions of the above theorem 
$\res{\ep_1}{\Ga_Z}\equiv\res{\ep_2}{\Ga_Z}\not\equiv 1$ but $v\neq 1$,
then  $(\Ga_1,\ep_1)$ and $(\Ga_2,\ep_2)$ are ``representation equivalent
up to multiplicity $v$'', i.e. 
$L^2_{\ep_1}(\Ga_1 \bs G)$ and 
$L^2_{\ep_2}(\Ga_2 \bs G)$ have the same irreducible components
and for every irreducible component $H$ with multiplicities $m_1,m_2$ we get 
$v=m_1 /m_2$.
\end{remark}

Let us compare this result to Theorem \ref{heisenspec}.
Using the notation of Section 3 for the 3-dimensional Heisenberg groups
write $M(r_i,d,T)= \Ga_i \bs G$. Then $v=r_1 /r_2$.
The remark explains why the eigenvalues of the Dirac operator
in Theorem \ref{heisenspec}~B do not depend on $r$ and 
why the multiplicity is proportional
to $r$.    

\proof{of Theorem \ref{nonsingquot} and Remark \ref{multremark}}
If $(\Ga_1\bs G, \ep_1)$ and $(\Ga_2\bs G, \ep_2)$ are representation 
equivalent, then $L^2_{\ep_1}(\Ga_1 \bs G)$ and $L^2_{\ep_2}(\Ga_2 \bs G)$
are also equivalent as $\Ga_Z$-modules. 
Since the action of $g_0 \in \Ga_Z$ on 
$L^2_{\ep_i}(\Ga_i \bs G)$ is just multiplication by $\ep_i(g_0)$ we get 
  $$\res{\ep_1}{\Ga_Z}\equiv\res{\ep_2}{\Ga_Z}.$$
Let $\Tau$ be the set of coadjoint orbits of $\gg^*$. Via Kirillov theory 
$\Tau$ parametrizes the set of irreducible unitary representations of $G$.
We write elements of $\Tau$ in the form $[\tau]$ with $\tau \in \gg^*$.
The irreducible $G$-module corresponding to $[\tau] \in \Tau$ will be denoted by 
$H_{[\tau]}$. 
We write $L^2_{\ep_i}(\Ga_i\bs G)$ as the discrete sum of irreducible 
representations 
  $$L^2_{\ep_i}(\Ga_i\bs G)=\bigoplus_{[\tau] 
    \in \Tau}m_i([\tau]) H_{[\tau]}$$
where $m_i([\tau])\in \NN_0$ denotes the multiplicity of $H_{[\tau]}$.
Let $\zz$ be the Lie algebra of the center $Z(G)$.
The action of $Z(G)$ on $H_{[\tau]}$ is trivial 
if and only if $\res{\tau}{\zz}\equiv 0$. 
Therefore
  $$\Hp i:=\bigoplus_{[\tau]\in \Tau \atop \res\tau\zz\equiv 0}
    m_i([\tau])H_{[\tau]} \qquad i=1,2 $$ 
is the space of $L^2$-functions acted trivially upon by $Z(G)$.
The orthogonal complement is
 $$H^\perp_i:=\bigoplus_{[\tau]\in \Tau \atop \res\tau\zz\not\equiv 0}
    m_i([\tau])H_{[\tau]} \qquad i=1,2\, .$$ 

The $G$-modules $L^2_{\ep_i}(\Ga_i\bs G)$ for $i=1,2$ 
are isomorphic if and only if $\Hp 1$ is isomorphic to $\Hp 2$ and 
$H^{\perp}_1$ is isomorphic to $H^{\perp}_2$.
In a first step we investigate under which conditions $\Hp 1$ and
$\Hp 2$ are isomorphic and then we treat $H^{\perp}_1$ and $H^{\perp}_2$. 

In the case $\res{\ep_i}{\Ga_Z}\not\equiv 1$ we have  $\Hp i=\{0\}$. 
Therefore $\Hp 1$ and $\Hp 2$ are trivially 
isomorphic $G$-modules.

In the case $\res{\ep_i}{\Ga_Z}\equiv 1$ elements in 
$\Hp i$ are exactly those functions in $L^2(\Ga_i\bs G)$ which
are pullbacks of functions in $L^2_{\ol\ep_i}(\ol\Ga_i\bs \ol G)$.
Thus we get a module-isomorphism (for both groups $G$ and $\ol G$)
  $$\Hp i \stackrel{\sim}{\longrightarrow} L^2_{\ol\ep_i}(\ol\Ga_i\bs \ol G).$$
Therefore in this case $\Hp 1$ and $\Hp 2$ are isomorphic
$G$-modules if and only if $(\ol\Ga_1,\ol\ep_1)$ and  $(\ol\Ga_2,\ol\ep_2)$
are representation equivalent for $\ol G$.

It remains to study the $H^\perp_i$. 
In Lemma~\ref{moorelemma} below we will show 
that $m_i([\tau])=c_{G,\tau} \cdot\vol(\ol\Ga_i\bs \ol G)$ 
for any $[\tau]\in \Tau$
with $\res\tau\zz\not \equiv 0$  where $c_{G,\tau}$
only depends on $G$, on $\tau$, and on the volume form on $\ol G$ but
not on the lattice $\Ga_i$.
This will complete the proof of the theorem:
Suppose that condition (1) or (2) in the theorem is fulfilled.
Because of Corollary \ref{volcor} we then know that the volume quotient $v$ 
is 1, and therefore $H^\perp_1$ and $H^\perp_2$ are isomorphic $G$-modules.
Using the statement for $v\neq 1$ we get Remark \ref{multremark}.
\qed

To complete the proof of Theorem~\ref{nonsingquot} and 
Remark~\ref{multremark} it remains to show Lemma~\ref{moorelemma}.
The lemma is a consequence of \cite{mw}. 
In order to simplify notation we drop the index $i$ 
in the formulation from now on.

For any $\rho\in \zz^*$ we define a bilinear form on $\ol{\gg}$ by
  $$ b_\rho (\ol X, \ol Y):= \rho ( [X,Y]),$$
where $X$ and $Y$ are vectors in $\gg$ and $\ol X$ and $\ol Y$ 
their images in $\ol{\gg}$.
As $\gg$ is strictly nonsingular, $b_\rho$ is nondegenerate for all 
$\rho\neq 0$ \cite{mw,g1}. 
Therefore $\dim \ol{\gg}$ is even and we suppose $\dim \ol{\gg}=2d$. 
We define the affine lattice $L_{\Ga,\ep}$ by 
  $$ L_{\Ga,\ep}:= \left\{ \rho \in \zz^* \, |\, 
     e^{2\pi i \rho(z)} = \ep(z) \quad
     \forall z \in \log(\Ga_Z)\right\}$$

\begin{lemma}\label{moorelemma}
If $\rho:=\res{\tau}{\zz} \in L_{\Ga,\ep}\ohne \{0\}$ then
  $$m([\tau])=\Big| \int\limits_{\ol\Ga\bs \ol G}
    \underbrace{b_\rho \wedge \ldots \wedge b_\rho}_{\mbox{$d$-times}}\Big|
    = \vol (\ol\Ga \bs \ol G)\cdot c_{G,\tau}$$
where $c_{G,\tau}$ only depends on $G$, on $\tau$ and on the volume
form on $\ol G$.

If $\rho:=\res{\tau}{\zz} \in \zz^*\ohne (L_{\Ga,\ep}\cup \{0\})$
then $m([\tau])=0$.
\end{lemma}

The lemma does not make any statement for $\rho=0$.
The reader may notice that for strictly nonsingular Lie groups the orbit
of $\res{\tau}{\zz}\neq 0$ is the set of all $\tau'\in \gg^*$ with 
$\res{\tau}{\zz}\equiv\res{\tau'}{\zz}$.

\proof{of Lemma \ref{moorelemma}} In the case $\ep\equiv 1$ the lemma 
is just a consequence of \cite[Theorem 7]{mw}. 
In the case $\ep\not\equiv 1$  we set $\Ga':=\ker(\ep)$ and we decompose
  $$ L^2_1(\Ga' \bs G) = L^2_1(\Ga \bs G ) \oplus L^2_\ep(\Ga\bs G)$$
as an orthogonal sum of $G$-modules.
We already now the lemma for $(\Ga,1)$ and $(\Ga',1)$. 
Using the fact that
every factor occurs only finitely many times,  we can calculate the 
multiplicity in  $L^2_\ep(\Ga\bs G)$ as the difference of the multiplicities
of $ L^2_1(\Ga' \bs G)$  and  $L^2_1(\Ga \bs G )$.
If $\ep$ is projectable, then 
$L_{\Ga,\ep}=L_{\Ga,1}=L_{\Ga',1}$, but 
$\vol(\ol\Ga'\bs \ol G)=2\vol(\ol\Ga\bs \ol G)$ 
and if $\ep$ is nonprojectable then 
$L_{\Ga',1}=L_{\Ga,1} \dot{\cup} L_{\Ga,\ep}$ 
but $\vol(\ol\Ga\bs \ol G)= \vol(\ol\Ga'\bs \ol G)$.
This proves the lemma for general $\ep$.
\qed

\begin{theorem}
There is a family of nilmanifolds having 

1. constant spectrum of the Laplacian on functions, 

2. nonconstant spectrum of the Laplacian on 1-forms

3. constant marked length spectrum

4. nonconstant spectrum of the Dirac operator for projectable 
spin structures,

5. constant spectrum of the Dirac operator for nonprojectable
spin structures.
\end{theorem}

\proof{}
Example I and II of \cite{g2} will have all these properties 
for suitably chosen lattices.
Properties 1--3 have been shown in \cite{g2}. We will prove 4 and 5
for Example I. The arguments for Example II are analogous.
To describe Example I let $\gg$ be the $7$-dimensional Lie algebra with
orthonormal basis 
$\{X_1,X_2,X_3,X_4,Z_1,Z_2,\cZ\}$ and Lie bracket

\begingroup
\catcode`@=\active
\def@{$ & $ = $ & $}
\begin{center}
\begin{tabular}{rccclcl}
$[X_1,X_2] @ [X_3,X_4] @ Z_1+ \cZ $\\
$[X_1,X_3] @ [X_4,X_2] @ Z_2      $\\
$[X_2,X_3] @ [X_1,Z_1] @ [Z_2,X_4] @ \cZ.$\\
\end{tabular}
\end{center}
\endgroup

All other Lie brackets of basis vectors shall be zero.
This Lie algebra is strongly nonsingular and $\cZ$ generates 
the center of $\gg$.
The matrix
  $$
\pmatrix{ 
\cos s & 0      & 0       & \sin s& 0        & 0      & 0\cr
0      & \cos 2s& -\sin 2s& 0     & 0        & 0      & 0\cr
0      & \sin 2s& \cos 2s & 0     & 0        & 0      & 0\cr
-\sin s& 0      & 0       & \cos s& 0        & 0      & 0\cr
0      & 0      & 0       & 0     & \cos s   & -\sin s& 0\cr
0      & 0      & 0       & 0     & \sin s   & \cos s & 0\cr
0      & 0      & 0       & 0     & -1+\cos s& -\sin s& 1\cr
}
  $$
defines a family of Lie algebra automorphisms $\ph_s$ of $\gg$ 
and induces a family of Lie group automorphisms $\Ph_s$ of $G$, 
i.e.\ ${\Ph_s}_*=\ph_s$.

The center $Z(G)$ is fixed under the automorphism $\Ph_s$.
Therefore $\Ph_s$ induces 
an automorphism $\ol\Ph_s$ of $\ol G = G / Z(G)$ with commuting diagram
\vspace{-1cm}
$$
\renewcommand{\arraystretch}{5.5}
\hfill
\begin{array}{r@{\hskip 2cm}l}
\rnode{a}{G} & \rnode{b}{G} \\
\rnode{c}{\ol G}&\rnode{d}{\ol G .}  \\
\end{array}
\psset{nodesep=3pt}
\ncline{->}{a}{b}\Aput{\Ph_s}
\ncline{->>}{a}{c}\Aput{\pi}
\ncline{->>}{b}{d}\Aput{\pi}
\ncline{->}{c}{d}\Aput{\ol\Ph_s}
\hfill
$$
These $\ol\Ph_s$ also act isometrically.

Let $\Ga$ be an arbitrary discrete cocompact subgroup of $G$. 
We write
$\Ga_s:=\Ph_s(\Ga)$, $\ol\Ga=\pi(\Ga)$ 
and $\ol\Ga_s:=\ol\Ph_s(\ol\Ga)$.

We now prove that the family 
$M_s:= \Ga_s\bs G$, $s \in \RR$, satisfies the theorem.
Gornet showed in \cite{g2} that this family is isospectral 
for the Laplacian on functions, but not for the Laplacian on 1-forms.
Furthermore she showed 
that $\Ph_s$ marks an isomorphism of the length spectrum
from $M_0$ to $M_s$.

We have to show 4 and 5 and we start with 5. 
Let $\ep : \Ga \to \{-1,1\}$ describe a nonprojectable spin structure on $M_0$.
As we want the spin structure to be chosen continuously in $s$, 
the spin structure on $M_s$ is given by $\ep \circ {\Ph_s}^{-1}$, 
and therefore also nonprojectable.
Since $\ol\Ph_s$ is a family of isometric automorphisms, 
the volume of $\ol\Ga_s \bs \ol G$ is independent of $s$.
Hence the second condition in Theorem \ref{nonsingquot} is satisfied
and therefore $\left(\Ph_s \circ \Ga, \ep \circ {\Ph_s}^{-1}\right)$ 
is a pairwise representation equivalent family. 
Theorem \ref{repequi} implies that the spectrum of the Dirac operator 
does not depend on $s$.

To show 4 let $M_0$ carry a projectable spin structure described by
$\ep:\Ga \to \{-1, 1\}$.
Let $\ol\ep:\ol\Ga \to \{-1,1\}$ describe
the induced spin structure on $\ol\Ga \bs \ol G$. 
By continuity as above 
$\ep_s:=\ep \circ {\Ph_s}^{-1}$ and 
$\ol\ep_s:=\ol\ep \circ {\ol\Ph_s}^{-1}$ describe 
the spin structures on $M_s$ and $\ol M_s:= \ol\Ga_s \bs \ol G$ resp.
 
Just as in the proof of Theorem \ref{repequi} we write the space of 
$L^2$-spinors on $M_s$ as 
  $$L^2\Big(\Si_{\ep_s}(\Ga_s\bs G)\Big)= 
    L^2_{\ep_s} \left(\Ga_s \bs G\right) \otimes_\CC \Si_7.$$
We can assume for our calculations that 
$$
\om := \ga (E_1)\cdot  \ga (E_2)\cdot\ldots\cdot \ga (E_7)
$$ 
acts as $-\Id$ on $\Si_7$. 
As in the proof of Theorem \ref{nonsingquot} we write 
  $$L^2_{\ep_s} \left( \Ga_s \bs G\right) = 
    \bigoplus_{[\tau] \in \Tau}m_s([\tau]) H_{[\tau]}.$$
The Leibniz rule \ref{leibniz} shows that 
$H_{[\tau]}\otimes \Si_7$ is a subspace invariant under 
$\na_{E_j}: L^2 \left( \Si_{\ep_s}( \Ga _s \bs G)\right)\to  
L^2 \left( \Si_{\ep_s}( \Ga _s \bs G)\right)$ and therefore
it is also invariant under the Dirac operator 
$D= \sum_{j=1}^7(\id \otimes \ga(E_j)) \cdot \na_{E_j}$.
Hence
  $$ \spec D^{M_s} = \bigcup_{[\tau] \in \Tau} 
     \spec \res{D}{H_{[\tau]}\otimes \Si_7}$$
and the multiplicity of each eigenvalue $\la$ of $D^{M_s}$ is the sum of 
$m_s([\tau])$ times the multiplicity of the eigenvalue $\la$ in 
$\res{D}{H_{[\tau]}\otimes \Si_7}$.

To prove that $\spec D^{M_s}$ is not constant in $s$ it 
is sufficient to show that there is one eigenvalue $\la_s$ in 
$\spec D^{M_s}$ which is continuous and nonconstant in $s$.

Therefore we will show that there is a 
$\tau \in \gg^*$ vanishing on $[\gg,\gg]$ such that 
$\tau_s:=\tau\circ {\Ph_s}^{-1}$ satisfies
$m_s([\tau_s])>0$ and $\res{D}{H_{[\tau_s]}\otimes \Si_7}$ has at least
one nonconstant eigenvalue depending continuously on $s$.

Because of $\res{\tau_s}{[\gg,\gg]}\equiv 0$ we know that $H_{[\tau_s]}$ 
is one-dimensional and $\exp X\in G$ acts on $H_{[\tau_s]}$ by multiplication
with $e^{2\pi i \tau_s(X)}$.
Elements $f$ of $H_{[\tau_s]}$ satisfy
$f(\exp X)=c \cdot e^{2\pi i \tau_s(X)}$ for all $X \in \gg$ with $c\in \CC$. 
The multiplicity $m_s([\tau_s])$ is $1$
if $e^{2\pi i \tau_s(X)}= \ep_s(\exp X)$ for any $\exp X \in \Ga_s$ 
and otherwise the multiplicity is $0$. This condition is independent of $s$.

We will compute the determinant of 
$\res{D}{H_{[\tau_s]}\otimes \Si_7}$. 
Writing $E_1:=X_1$, \dots, $E_4:=X_4$, $E_5:=Z_1$, $E_6:=Z_2$ and $E_7:=\cZ$,
we apply (\ref{nablaspinformel}) and get for $S \in \Si_7$ and 
$f \in H_{[\tau_s]}$:
\begin{eqnarray*}
  D^{M_s}(fS)& = & \sum_i \left( R_*(E_i) f\right) \ga(E_i) S
  + {1\over 4} f \sum_{i,j,k} \Ga_{ij}^k \ga(E_i)\ga(E_j)\ga(E_k)S \cr
& = & \sum_{i=1}^4 2 \pi i \tau_s(X_i) f \ga(E_i) S + f A(S) \cr
\end{eqnarray*}
where 
$A:=(1/4) \sum_{i,j,k} \Ga_{ij}^k \ga(E_i)\ga(E_j)\ga(E_k) 
\in \End(\Si_7)$
is constant in $s$. 
Using the Koszul formula one can express $A$ in a suitable basis by 
the matrix

$$ 
\begingroup
\catcode`+=\active
\def+{\phantom{-}}
A={1\over 4}
\pmatrix{
-2 & +0 & +0 & +1 & +2 & +{i} & +{i} & -2 \cr
+0 & +0 & +1 & +0 &  -{i} & +0 & +0 &  -{i} \cr
+0 & +1 & +0 & +0 & +{i} & +0 & +0 & +{i} \cr
+1 & +0 & +0 & +2 & +2 &  - {i} &  - {i} & -2 \cr
+2 & +{i} &  - {i} & +2 & +2 & +0 & +0 & -1 \cr
-{i} & +0 & +0 & +{i} & +0 & +0 & -1 & +0 \cr
-{i} & +0 & +0 & +{i} & +0 & -1 & +0 & +0 \cr
-2 & +{i} &  - {i} & -2 & -1 & +0 & +0 & -2
}
\endgroup$$
Then 
\begin{eqnarray*}
\det \res{D}{H_{[\tau_s]}\otimes \Si_7}
& = & \det\bigg[ A + 2 \pi i \sum_{i=1}^4 \tau_s(X_i)\ga(E_i) \bigg]
\end{eqnarray*}  
and 
\begin{eqnarray*}
  \de_{\tau}& := & 
     \stelle{d\over ds}{s=0}\det \res{D}{H_{[\tau_s]}\otimes \Si_7}\cr
& = &  - \, \frac {1}{64}{ \pi}^{2}\, \left( 
\! \, - {\tau(X_1)}\,{\tau(X_2)} + {\tau(X_3)}\,{\tau(X_4)}\, 
\!  \right)\\  
&& \Big( 4096\,{ \pi}^{4}\,{\tau(X_1)}^{4} 
+ 4096\,{ \pi}^{4}\,{\tau(X_2)}^{4}
+ 4096\,{ \pi}^{4}\,{\tau(X_3)}^{4} 
+ 4096\,{ \pi}^{4}\,{\tau(X_4)}^{4}\\
&& \mbox{}+ 8192\,{ \pi}^{4}\,{\tau(X_1)}^{2}\,{\tau(X_2)}^{2}
+ 8192\,{ \pi}^{4}\,{\tau(X_1)}^{2}\,{\tau(X_3)}^{2} 
+ 8192\,{ \pi}^{4}\,{\tau(X_2)}^{2}\,{\tau(X_3)}^{2}\\
&& \mbox{}+ 8192\,{ \pi}^{4}\,{\tau(X_1)}^{2}\,{\tau(X_4)}^{2}
+ 8192\,{ \pi}^{4}\,{\tau(X_2)}^{2}\,{\tau(X_4)}^{2} 
+ 8192\,{ \pi}^{4}\,{\tau(X_3)}^{2}\,{\tau(X_4)}^{2}\\
&& \mbox{}
+ 512\,{ \pi}^{2}\,{\tau(X_3)}\,{\tau(X_1)}  
+ 128\,{ \pi}^{2}\,{\tau(X_4)}^{2} + 512\,{ \pi}^{2}\,
{\tau(X_4)}\,{\tau(X_2)} + 128\,{ \pi}^{2}\,{\tau(X_2)}^{2} \\
&& \mbox{}  + 128\,{ \pi}^{2}\,{\tau(X_1)}^{2}
+ 128\,{ \pi}^{2}\,{\tau(X_3)}^{2} 
 + 1 \Big). 
\end{eqnarray*}
The set 
  $$\Big\{\big(\tau(X_1),\tau(X_2),\tau(X_3),\tau(X_4)\big)\,|\,
    \tau\in \gg^*, \res{\tau}{[\gg,\gg]}\equiv 0, m([\tau])>0\Big\}$$
is an affine lattice in $\RR^4$ 
(i.e.\ up to translation a discrete cocompact subgroup). 
But the only polynomial vanishing on an affine lattice of $\RR^4$ is the 
zero polynomial. 
Therefore there exist $\tau \in \gg^*$ with 
$\res{\tau}{[\gg,\gg]}\equiv 0$, $m([\tau])>0$ and $\de_\tau\neq 0$. 
\qed
 

\vspace{1cm}

\noindent
Authors' address:

\vspace{0.3cm}

\noindent
Mathematisches Institut\\
Universit\"at Freiburg\\
Eckerstr.~1\\
79104 Freiburg\\
Germany

\vspace{0.3cm}

\noindent
e-mail:
\newline 
{\typewriter ammann@mathematik.uni-freiburg.de}
\newline
\noindent 
{\typewriter baer@mathematik.uni-freiburg.de}

\end{document}